\documentclass[12pt]{article}
\usepackage{hyperref,amsmath,amssymb}
\usepackage{eso-pic}
\usepackage{makeidx}
\usepackage{subfiles}

\makeatletter
\def\printnotation{{%
\def\indexname{Index of notation}
\begin{theindex}
\@input{\jobname.ntn}
\end{theindex}
}}
\makeatother

\makeglossary
\usepackage{amsfonts,amssymb,amscd,amsmath,latexsym,amsbsy, setspace}
\usepackage{amsfonts,amssymb,graphicx,amsthm,mathtools}
\usepackage{mathrsfs}
\usepackage{subfiles}
\usepackage[all,cmtip,curve, knot,frame]{xy} 
\usepackage{amssymb}
\usepackage{amsfonts}
\usepackage{latexsym}
\usepackage{epstopdf}
\usepackage{tikz}
\usepgflibrary{shapes.geometric}
\usetikzlibrary{matrix}

\usepackage{enumitem}
\setlist{noitemsep}

\newtheorem{theorem}{Theorem}[section]
\newtheorem{lemma}[theorem]{Lemma}

\newtheorem{proposition}[theorem]{Proposition}
\newtheorem{corollary}[theorem]{Corollary}
\theoremstyle{definition}
\newtheorem{definition-proposition}[theorem]{Definition-Proposition}
\newtheorem{definition}[theorem]{Definition}
\newtheorem{example}[theorem]{Example}
\newtheorem{remark}[theorem]{Remark}

\newcommand{\LFP}{\operatorname{\mathbf{Pr}}}
\newcommand{\LFPd}{\operatorname{\mathbf{Pr}}^\circ}

\newcommand{\colim}{\operatorname{colim}}

\newcommand{\mop}{{mop}}
\newcommand{\bop}{{bop}}

\newcommand{\Vect}{\operatorname{Vect}}

\newcommand{\End}{\operatorname{End}}
\newcommand{\Hom}{\operatorname{Hom}}
\newcommand{\IHom}{\underline{\operatorname{Hom}}}

\newcommand{\Rep}{\operatorname{Rep\,}}

\newcommand{\K}{\mathbf{k}}

\newcommand{\cA}{\mathcal{A}}
\newcommand{\cX}{\mathcal{X}}
\newcommand{\cB}{\mathcal{B}}
\newcommand{\cC}{\mathcal{C}}

\newcommand{\cD}{\mathcal{D}}
\newcommand{\cF}{\mathcal{F}}

\newcommand{\cM}{\mathcal{M}}

\newcommand{\cN}{\mathcal{N}}
\newcommand{\cS}{\mathcal{S}}
\newcommand{\D}{\mathbb{D}}
\newcommand{\into}{\hookrightarrow}

\newcommand{\Lin}{\operatorname{Lin}}
\newcommand{\Rex}{\operatorname{\mathbf{Rex}}}

\newcommand{\Fun}{\operatorname{\operatorname{Fun}}}

\newcommand{\cE}{\mathcal{E}}

\newcommand{\ot}{\otimes}
\newcommand{\bt}{\boxtimes}
\newcommand{\id}{\operatorname{id}}

\newcommand{\RR}{\mathbb{R}}

\newcommand{\coev}{\operatorname{coev}}
\newcommand{\ev}{\operatorname{ev}}

\newcommand{\oo}{\infty}

\newcommand{\Tens}{\operatorname{\mathbf{Tens}}}
\newcommand{\BrTens}{\operatorname{\mathbf{BrTens}}}
\newcommand{\RigidTens}{\operatorname{\mathbf{RigidTens}}}
\newcommand{\RigidBrTens}{\operatorname{\mathbf{RigidBrTens}}}
\newcommand{\BrFus}{\operatorname{\mathbf{BrFus}}}
\newcommand{\Fus}{\operatorname{\mathbf{Fus}}}

\newcommand{\coker}{\operatorname{ker}}

\newcommand{\uch}{\underline{\operatorname{Ch}}}
\newcommand{\ch}{\operatorname{Ch}}

\newcommand{\modu}{\operatorname{-mod}}
\newcommand{\modr}{\operatorname{mod-}\!}
\newcommand{\modul}{\operatorname{mod}}

\newcommand{\un}{\mathbf{1}}
\usepackage{soul}

\usepackage{textcomp}

\def\HH{\hbox{${\mathcal H}$\kern-5.2pt${\mathcal H}$}}

\begin{document}
\title{On dualizability of braided tensor categories}
\author{Adrien Brochier, David Jordan, Noah Snyder}

\maketitle
\begin{abstract}
We study the question of dualizability in higher Morita categories of locally presentable tensor categories and braided tensor categories.  Our main results are that the 3-category of rigid tensor categories with enough compact projectives is 2-dualizable, that the 4-category of rigid braided tensor categories with enough compact projectives is 3-dualizable, and that (in characteristic zero) the 4-category of braided fusion categories is 4-dualizable. Via the cobordism hypothesis, this produces respectively 2, 3 and 4-dimensional framed local topological field theories.  In particular, we produce a framed 3-dimensional local TFT attached to the category of representations of a quantum group at any value of $q$.
\end{abstract}
\setcounter{tocdepth}{2}
\tableofcontents

\section{Introduction}
This paper establishes dualizability results for the higher Morita categories of locally presentable tensor categories and braided tensor categories.  Roughly speaking, a higher category is called \emph{fully dualizable} if its objects have duals and its morphisms have adjoints at all levels.  More generally, a higher category is $k$-dualizable if all objects have duals and all morphisms of degree less than $k$ have adjoints.  The \emph{cobordism hypothesis} of Baez-Dolan~\cite{Baez1995} hypothesizes that local (or fully extended) topological field theories in a given dimension, valued in some target higher symmetric monoidal category, are classified by its fully dualizable objects, i.e. those objects lying in some fully dualizable subcategory.  Here local means that it is extended all the way down to points, i.e. the TFTs assign objects of the target to 0-manifolds, 1-morphisms to 1-dimensional bordisms, 2-morphisms to 2-dimensional bordisms with corners, etc.  Lurie's influential work~\cite{Lurie2009} outlined a proof of the cobordism hypothesis; several more recent works~\cite{Ayala2017a,Ayala2015a} have been aimed at giving an independent proof.  In light of this classification, it is therefore a very interesting question to produce examples of fully dualizable higher categories, and $k$-dualizable higher categories more generally.

Of particular interest are the \emph{higher Morita theories of $E_n$-algebras}.  These generalize the classical Morita theory -- the 2-category of associative algebras, bimodules, and bimodule homomorphisms -- to the setting of $E_n$-algebras in some symmetric monoidal $m$-category $\cS$, i.e. algebras in $\cS$ over the $n$-little discs operad.  An $E_n$-algebra may be regarded as carrying $n$ mutually distributing associative algebra structures in $\cS$, and so the resulting Morita theory of $E_n$-algebras is most naturally regarded as an $(n+m)$-category\footnote{Throughout the paper, by ``$n$-category" we will mean more precisely an $(\infty,n)$-category, in the complete Segal space axiomatization.  We will say ``discrete $n$-category" when the only invertible higher morphisms are identites.  However, following the remark in Section I.5 of \cite{Douglas2013} all our main arguments take place in discrete $2$-categories and so are model independent.}.  The first general construction of a higher Morita theory was the $(n+1)$-category of $E_n$-algebras in a (possibly non-discrete) $1$-category $\cS$.  This construction was outlined in~\cite{Lurie2009}, and carried out in independent works \cite{Haugseng2017,Scheimbauer2014}.  The further construction of the $(n+m)$-categorical higher Morita theory of $E_n$-algebras in an $m$-category $\cS$ was given in \cite{Johnson-Freyd2017}.  

When we specialize the coefficients $\cS$ appearing in these constructions to be the (discrete) 2-category $\LFP$ of locally presentable $\K$-linear categories, we obtain a (discrete) $3$-category $\Tens$ of tensor categories, regarded as $E_1$-algebras in $\LFP$, and a (discrete) $4$-category $\BrTens$ of braided tensor categories, regarded as $E_2$-algebras in $\LFP$.  Our main results in this paper consist of identifying natural subcategories of $\Tens$ and $\BrTens$ consisting of \emph{rigid} tensor and braided tensor categories, respectively, with enough compact projectives, and proving that these are $2$- and $3$-dualizable, respectively.  In the braided tensor case, we further identify a subcategory of \emph{braided fusion categories}, and we show this is fully dualizable.  Taken together with the cobordism hypothesis, our results give rise to new topological field theories in dimensions $2$, $3$, and $4$.  In the tensor case this generalizes results of \cite{Douglas2013} for \emph{finite} rigid tensor categories, is a non-derived analogue of \cite{Ben-Zvi2009,Gaitsgory2015}, and extends to tensor and braided tensor categories the dualizability results of \cite{Brandenburg2015} for presentable categories.  

The rest of the introduction is outlined as follows.  We begin by recalling what the Morita categories of tensor and braided tensor categories are, briefly discuss the notion of rigidity, as it appears in the infinite/presentable setting, and state our main results. Then we discuss the relationship of our results to previous results on dualizability of tensor categories.  We next briefly discuss what we expect to be true about oriented TFTs attached to braided tensor categories, which we will return to in future work. Finally we discuss the relationship of the $3$-dimensional TFTs coming from braided tensor categories via the cobordism hypothesis to a number of well-known constructions, some of which are at present only conjecturally well-defined.

\subsection{Morita categories} Our primary objects of interest are tensor categories and braided tensor categories, and their higher Morita theories.  Following the general framework laid out in \cite{Haugseng2017, Scheimbauer2014, Johnson-Freyd2017}, Section \ref{sec:BrTens} of the present paper is devoted to the case $\cS=\Pr$, where we can spell out explicitly the data of all higher morphisms, and the composition laws.

\begin{definition-proposition}\label{def:tens} There exists of a $3$-category $\Tens$, whose:
\begin{itemize}
\item objects are tensor categories.
\item 1-morphisms are bimodule categories between such.
\item 2-morphisms are bimodule functors between such.
\item 3-morphisms are natural transformations between such.
\end{itemize}
\end{definition-proposition}

The symmetric monoidal product on $\Tens$ is the so-called Deligne--Kelly tensor product of the underlying categories, equipped with its natural tensor structure.  The composition of $1$-morphisms is the balanced (or relative) Deligne--Kelly tensor products of bimodules.

Compositions of $2$- and $3$-morphisms are composition of functors and of natural transformations in the usual sense.

Similarly, we have:

\begin{definition-proposition}\label{def:brtens} There exists a $4$-category $\BrTens$, whose:
\begin{itemize}
\item objects are braided tensor categories.
\item 1-morphisms are tensor categories with central structures.
\item 2-morphisms are centered bimodule categories.
\item 3-morphisms are bimodule functors of such.
\item 4-morphisms are natural transformations of such.
\end{itemize}
\end{definition-proposition}

The symmetric monoidal product in $\BrTens$ is again the Deligne--Kelly tensor product of the underlying categories.  The composition of $1$- and $2$-morphisms is defined using the balanced Deligne-Kelly tensor product.  Composition of $3$- and $4$-morphisms is composition of functors and natural transformations in the usual sense.  See Section \ref{sec:BrTens} for complete definitions.

\subsection{Rigidity} With these definitions in hand, our main results identify certain appropriately dualizable sub-categories of $\Tens$ and $\BrTens$. 
 A key notion in these theorems is that of rigidity.  In the setting of finite tensor categories, rigidity of a tensor category is the requirement that all its objects have left and right dual objects, equipped with evaluation and coevaluation maps mimicking those associated to the dual of a vector space.   In the infinite setting, some care is needed to define the correct notion of rigidity: for instance we wish to include the category of vector spaces as an example of a rigid tensor category, although it contains infinite-dimensional vector spaces, which are not dualizable.  This leads us to: 

\begin{definition-proposition}
Suppose that a tensor category $\cA$ has enough compact projectives.  Then the following conditions on $\cA$ are equivalent:
\begin{enumerate}
\item All compact projective objects of $\cA$ are left and right dualizable.
\item A generating collection of compact projective objects of $\cA$ are left and right dualizable.
\item The multiplication functor, $T: \cA\bt\cA\to\cA$ has a co-continuous right adjoint $T^R$, and the canonical lax bimodule structure on $T^R$ is strong.
\end{enumerate}
 We will say that $\cA$ is \emph{cp-rigid} if it has enough compact projectives, and if any of the above conditions is satisfied.
\end{definition-proposition}

\begin{remark} Note that we only define cp-rigidity of $\cA$ under the standing assumption that the underlying category has enough compact projectives.  We will follow this convention throughout the paper, so that a tensor or braided tensor category declared to be cp-rigid in particular has enough compact projectives.  Nonetheless in the main theorems we explicitly include the enough compact projectives condition in order to avoid any possible confusion.
\end{remark}

\begin{remark} The reader should compare the notion of cp-rigidity above with that of~\cite[Appendix D]{Gaitsgory2015}, in the setting of dg-categories.  There a dg-category is called rigid if it satisfies the third condition above and if the right adjoint to the unit is cocontinuous. Dropping this last condition leads to what is called ``semi-rigid'' in~\cite{Ben-Zvi2009}.
\end{remark}

\begin{remark} \label{rem:Theo}
Another reasonable notion is compact-rigid, meaning that all compact objects have duals.  This more closely matches the traditional notion of rigidity appearing in finite tensor category literature, in the following sense: if a finite tensor category which is in the traditional sense (that all objects are dualizable), then its ind-completion is compact-rigid.  Clearly, compact-rigid implies cp-rigid, but the former is more restrictive.

Even in the finite setting cp-rigidity is better behaved in some important ways, as was pointed out to us by Theo Johnson-Freyd.  For example, if $A$ is a finite-dimensional, non-semisimple algebra in $\Vect$ then $A$-mod-$A$ is not compact-rigid even though it is Morita equivalent to $\Vect$.  By contrast, $A$-mod-$A$ is cp-rigid, because it is generated by the compact projective bimodule $A\otimes A$, which is dualizable (because it is projective -- and indeed free -- as both a left and right $A$-module).

The simplest instance to see the distinction is $A = \mathbb{C}[x]/x^2$; here it is not hard to see that the only indecomposable dualizable objects are the projective bimodule $A \otimes A$ and the trivial bimodule.  In particular, the bimodule $\mathbb{C}$ with both the left and right $x$ acting by $0$ is compact but not dualizable.
\end{remark}

\begin{remark}
Similarly cp-rigid is better behaved over imperfect fields than compact-rigid, since compact-rigid is not preserved by the Deligne-Kelly tensor product (See \cite[Prop. 5.17]{Deligne2007}) while cp-rigid is preserved.

A specific counterexample from \cite{Douglas2013} arises from an inseparable field extension $L/K$.  Recall that a finite dimensional module over finite dimensional algebra is dualizable if, and only if, it is projective; hence, $A$-mod-$A$ is compact-rigid if, and only if, every finite dimensional bimodule is projective as a left and as a right module which happens if, and only if, $A$ is semisimple.   Thus $L$-mod-$L$ is compact-rigid (since $L$ is semisimple) but its tensor square $L\text{-mod-}L \boxtimes L\text{-mod-}L \cong (L\otimes L)\text{-mod-}(L\otimes L)$ is not compact-rigid, since $L \otimes L$ is not semisimple. However, as pointed out to us by Theo Johnson-Freyd,  $(L\otimes L)\text{-mod-}(L \otimes L)$ is cp-rigid in the sense of this paper.
\end{remark}

\subsection{Main results}
We can now state our main results:
\begin{theorem}\label{thm:Tens-Informal} There exists a sub 3-category $\RigidTens\subset \Tens$, whose:
\begin{itemize}
	\item objects are \emph{cp-rigid} tensor categories with enough compact projectives
\item 1-morphisms are bimodule categories with enough compact projectives.
\item 2-morphisms are bimodule functors between such.
\item 3-morphisms are natural transformations between such.
\end{itemize}
Moreover, $\RigidTens$ is 2-dualizable.
\end{theorem}

\begin{theorem}\label{thm:BrTens-Informal} There exists a sub 4-category $\RigidBrTens \subset \BrTens$, whose:
\begin{itemize}
\item objects are \emph{cp-rigid} braided tensor categories with enough compact projectives.
\item $1$-morphisms are \emph{cp-rigid} tensor categories with central structures, and enough compact projectives,
\item $2$-morphisms are centered bimodule categories, with enough compact projectives,
\item $3$- and $4$-morphisms are as in $\BrTens$.
\end{itemize}
Moreover, $\RigidBrTens$ is 3-dualizable.
\end{theorem}

It is already known from \cite{Douglas2013} that fusion categories over a field of characteristic zero are $3$-dualizable.  We show that the analogous statement holds for $4$-dualizability in $\BrTens$.

\begin{theorem}  Over a field of characteristic zero, there exists a sub 4-category $\BrFus \subset \BrTens$, whose:
\begin{itemize}\label{thm:BrFus-Informal}
\item objects are braided fusion categories.
\item $1$-morphisms are fusion categories with compatible central structures.
\item $2$-morphisms are finite and semi-simple bimodule categories with compatible centered structures.
\item $3$-morphisms are \emph{compact preserving} functors between such.
\item 4-morphisms are natural transformations between such.
\end{itemize}
Moreover, $\BrFus$ is fully (i.e. 4-) dualizable.
\end{theorem}

As in \cite{Douglas2013} over a field of characteristic $p$, this theorem holds with semisimplicity replaced by the stronger condition of ``separability."

\begin{remark}
The statements of the main Theorems \ref{thm:Tens-Informal}, \ref{thm:BrTens-Informal}, and \ref{thm:BrFus-Informal} are somewhat simplified for the sake of exposition.  More precisely, we first prove in Section 4 that each desired subcategory is indeed closed under the monoidal structure and the composition of morphisms at each level.  We then construct adjoints in Section 5 to $k$-morphisms in the required degrees.
\end{remark}

\subsection{Summary of known dualizability results}
Now we recall what has been proved already regarding dualizability in higher Morita theories, and explain where our results fit.

Dualizability in low dimensions is a general and purely topological phenomenon.  It was proved in \cite{Lurie2009} that any $E_1$-algebra $\cA$ (in any $\cS$) is $1$-dualizable, with dual being the opposite algebra $\cA^{\mop}$.  (Here the $m$ stands for multiplication, as opposed to taking the opposite of composition or of the braiding.)  This result was generalized in \cite{Gwilliam2018} to arbitrary $n$, establishing that the entire pointed higher Morita category of $E_n$-algebras is fully $n$-dualizable; for objects this was proved in \cite{Scheimbauer2014}. This is also expected to hold for unpointed higher Morita categories of~\cite{Haugseng2017}.  In other words: it follows from purely topological considerations that every object of the Morita theory of $E_n$-algebras is $n$-dualizable, in particular that $\Tens$ is $1$-dualizable and $\BrTens$ is $2$-dualizable.

Turning attention to our ``coefficients" $\LFP$, it was proved in~\cite{Brandenburg2015} that $\cC\in\Pr$ is $1$-dualizable if it has enough compact projectives\footnote{The converse assertion, that this is a necessary condition, is proved there only for categories of comodules for a co-algebra; the status of the converse in general does not appear to be known.}.  It is well-known that $\cC\in\LFP$ is 2-dualizable if it is the ind-completion of a finite and semi-simple category (and the converse holds, under the assumption there are enough compact projectives, see \cite{Tillman1998} and the Appendix of \cite{Bartlett2015}).

In \cite{Douglas2013} it was shown that \emph{finite} tensor categories are $2$-dualizable and that in characteristic zero \emph{fusion} categories -- semi-simple finite rigid tensor categories -- are fully (i.e. $3$-) dualizable.  In characteristic $p$ it was shown that separable tensor categories (i.e. fusion categories of non-zero global dimension) are fully dualizable.  Here we show that the first of these results generalizes beyond the finite setting.\footnote{Although our results show finiteness is not necessary for the $2$-dualizability theorem in \cite{Douglas2013}, many of the other results and conjectures of \cite{Douglas2013} -- such as the topological description of the Radford isomorphism, the conjectured $SO(3)$ action on finite tensor categories, and the conjectured non-compact local TFT -- do use finiteness in an essential way.}
In the infinite setting we have no improved sufficient condition for $3$-dualizability beyond fusion (resp. separable) categories from \cite{Douglas2013}, see Remark \ref{rem:converse}.

A necessary condition for $\cA\in\Tens$ to be 2-dualizable, or for $\cA\in\BrTens$ to be $3$-dualizable, is that $\cA$ be $1$-dualizable in $\LFP$.  Absent a general characterization of dualizability in $\LFP$, we take the criterion of having enough compact projectives, from \cite{Brandenburg2015}, which is in any case a versatile one for proofs.  Our results assert that, having assumed this form of dualizability in $\LFP$, it suffices to further assume only that $\cA$ is cp-rigid, in order to obtain the next degree of dualizability.

\begin{table}\begin{center}\begin{tabular}{c|c|c|c}
 & $\LFP$ & $\Tens$ & $\BrTens$\\
\hline
1 & $\begin{array}{c}\textrm{Enough}\\\textrm{compact}\\ \textrm{projectives}\end{array}$ & Any & Any\\
\hline
2 & $\begin{array}{c}\textrm{Finite}\\ \textrm{semi-simple}\end{array}$& Cp-rigid & Any\\
\hline
3 & & Fusion & Cp-rigid\\
\hline
4 & & &Fusion\\
\end{tabular}\end{center}\caption{Sufficient conditions for dualizability.  Dualizability in low dimensions is topological, while dualizability in higher dimensions references dualizability conditions in $\LFP$.  Cp-rigidity provides the additional compatibility in middle dimensions.}\label{table}\end{table}

To summarize (see Table~\ref{table}), $n$-dualizability in higher Morita theories of $E_n$-algebras is guaranteed by very general theorems~\cite{Scheimbauer2014,Haugseng2017,Lurie2009}, which are really about the topological $E_n$-operads, and in particular apply for any coefficients $\cS$.  Meanwhile, \emph{full} dualizability in higher Morita theories imposes very strong restrictions -- in particular the fusion condition in the works~\cite{Douglas2013,Freed2012a,Freeda,Walker}.  Our results highlight rigidity as a sufficient condition for one additional level of dualizability in $\Tens$ and $\BrTens$.  The ability to drop the strong finiteness assumptions, while retaining one level of dualizability -- 2-dualizability of tensor categories, and 3-dualizability of braided tensor categories -- is important for applications to geometric representation theory and quantum algebra (see below), where categories of interest are typically neither finite nor semi-simple, and where invariants are already known for other reasons not to extend to top dimension.

\begin{remark}
It is also interesting to compare these results to dualizability for ordinary finite-dimensional algebras.  There $2$-dualizability is equivalent to separability, i.e. that the multiplication map has a splitting as a bimodule map.  In some sense separability (splitting of multiplication) looks quite similar to rigidity (cocontinuous right adjoint to multiplication), but we do not know a common generalization that includes both separability and rigidity.
\end{remark}

\begin{remark}\label{rem:converse}
It is natural to ask whether the sufficient conditions from Table \ref{table} can be improved to neccessary and sufficient conditions characterizing dualizability.  We make no claims in this direction in this paper, however we note that some of these conditions are certainly not necessary.  In particular, as pointed out to us by Dan Freed, if $A$ is a finite dimensional non-semisimple algebra then the tensor category $A$-mod-$A$ is not fusion, but is $3$-dualizable because it is Morita equivalent to $\Vect$ (see the proof of Theorem 7.12.11 \cite{Etingof2015}).
\end{remark}

\subsection{Applications to topological field theories}

Combining the corbordism hypothesis with Theorems \ref{thm:Tens-Informal}, \ref{thm:BrTens-Informal} and \ref{thm:BrTens-Informal} yields the following

\begin{corollary}\label{cor:TFTs} We have:

	\begin{enumerate}
		\item Every cp-rigid tensor category gives rise to a local, categorified, framed 2-dimensional TFT assigning that category to the framed point.
		\item Every cp-rigid braided tensor category gives rise to a local, categorified, framed 3-dimensional TFT assigning that category to the framed point.
		\item In the braided fusion case, the this framed 3-dimensional TFT extends to a framed 4-dimensional TFT.
	\end{enumerate}
\end{corollary}

Here ``categorified'' means that one obtains vector spaces, as opposed to numbers, for closed manifolds in top dimension, and these are equipped with actions of mapping class groups.

In the braided case, the resulting TFT produces categories for framed surfaces, which coincide (via the uniqueness assertion in the cobordism hypothesis) with the ``quantum character varieties'' computed in~\cite{Ben-Zvi2015} using the formalism of factorization homology.  Extending the construction of character varieties to dimension 3 was a primary motivation for this work.  The 3-dimensional part of the TFT assigns vector spaces to closed 3-manifolds (more generally, functors between the quantum character theories of incoming and outgoing boundary components, to manifolds with boundaries), and we expect those vector spaces, (resp. functors) to admit a description via skein modules.  Finally, in the braided fusion case the TFT produces numerical invariants for framed 4 manifolds, which we expect will be closely related to those obtained by Crane--Kauffmann--Yetter~\cite{Crane1997} in the same manner that the $3$-dimensional TFT from \cite{Douglas2013} is expected to relate to Turaev--Viro invariants \cite{Turaev1992, Barret1996, Turaev2010}.  In this section, we outline all of these expected connections and applications in more detail.

Firstly, let us remark that the second and third statements of Corollary \ref{cor:TFTs} (in particular, their oriented variant to be discussed below) have been anticipated for many years, at least under the assumption that the input braided tensor category is semi-simple.  That ribbon tensor categories in characteristic $0$ give local $3$-dimensional TFTs was known to Kevin Walker \cite{Walker} and essentially follows from \cite{Morrison2012,Morrison2011}, though these works use an alternative notion of higher category and of TFT, which are not easily translated to the more standard ones used in this paper.    A variation on their approach is outlined in \cite[\S9]{Johnson-Freyd2015}; 
his proposal uses an axiomatic framework compatible with ours, and seems a likely candidate to give an independent proof of Corollary \ref{cor:TFTs}.  We discuss this in more detail below.

The statement that ribbon fusion categories give a fully local $4$-dimensional TFT also appears in unpublished work of Freed--Teleman~\cite{Freed2012a,Freeda}, and in unpublished work of Walker~\cite{Walker}, so in particular it should follow from the cobordism hypothesis that they are $4$-dualizable.  Hence, in the braided fusion case, Theorem \ref{thm:BrFus-Informal} can be regarded as giving an independent proof of this widely expected result, but on the other side of the cobordism hypothesis correspondence from other arguments.

All three constructions are closely related, and should agree exactly in the braided fusion case, but there are some important differences in general.  Firstly, the kind of completion that Walker--Morrison use (``the finite representation category") does not agree with the free cocompletion in general and so differs from the TFTs in Johnson-Freyd's approach and ours.\footnote{The finite representation category $\Fun^{\mathrm{add}}(\cC^{\mathrm{op}}, \Vect^{fd})$ is neither the free cocompletion in the presentable setting $\Fun^{\mathrm{add}}(\cC^{\mathrm{op}}, \Vect)$, nor the free cocompletion in finitely cocomplete setting $\Fun^{\mathrm{add}}(\cC^{\mathrm{op}}, \Vect)^{\mathrm{compact}}$ (unless $\cC$ is finite).}  Secondly and more importantly, in the skein formalism proposed by \cite{Johnson-Freyd2015} (see Section \ref{sec:skeins} for a detailed discussion) one must restrict strand types to compact projective objects; in particular this includes the unit which corresponds to the empty strand, and thus by a standard argument implies that all dualizable objects are compact projective.  This restriction allows for semi-simple cp-rigid braided tensor categories (finite or infinite), such as quantum groups at generic parameters, or modularizations of quantum groups at roots of unity, but excludes non-semisimple finite braided tensor categories in the sense of \cite{Etingof2004}, as well as infinite non-semisimple braided tensor categories, such as the (non-modularized) quantum groups at roots of unity.  The reason the compact projective assumption appears in the skein theoretic approach is that, by construction, any category includes fully faithfully into its free co-completion as the sub-category of compact projective objects, so that any attempt to freely co-complete skein categories renders a category where the unit is compact projective.

\subsubsection{The $\mathrm{SO}(n)$ action and oriented TFTs} \label{sec:oriented}

By the cobordism hypothesis, the space of $n$-dualizable objects (i.e. the underlying groupoid or core of the category of the category of fully dualizable objects) is homotopy equivalent to the space of $n$-dimensional \emph{framed} TFTs.  By precomposing with changes of framing, the group $\mathrm{O}(n)$ acts on the space of TFTs, and thus on the core of the category of dualizable objects $\cC^{\times}$.  The second part of the cobordism hypothesis in Lurie's formulation says that the oriented field theories correspond precisely to \emph{homotopy fixed points} for the $\mathrm{SO}(n)$ part of this action.  An action of a topological group $G$ on the core a higher category $\cC$ is given by a map $BG\to B\mathrm{Aut}(\cC^{\times})$, and a homotopy fixed point consists of a trivialization of a component of this action; this involves both non-trivial conditions on the would-be fixed points, and the specification of data giving the null homotopy.

The applications we have in mind (discussed in the next three sections) require oriented theories, so in this section we briefly outline what we expect to be true about the $\mathrm{SO}(3)$-action and its homotopy fixed points for $\RigidBrTens$. We only sketch the main ideas here, as the details will appear in a forthcoming paper by us; these hinge in turn on work in progress by one of us and co-authors~\cite{Douglas}.  It is also interesting to ask the corresponding questions about the $\mathrm{SO}(2)$ action on $\RigidTens$, the $\mathrm{SO}(3)$ action on $\Fus$, and the $SO(4)$ action on $\BrFus$.  The first two cases are considered in~\cite{Douglas}, while for the final case we do not have an outline of a proof in mind at present, and so will only speculate below.

The $\mathrm{O}(n)$ action on $E_n$-algebras is well-understood in concrete terms.  Namely, $\mathrm{O}(n)$ acts on the $E_n$ operad itself by the standard action on the little discs, and the $\mathrm{O}(n)$ action on $E_n$-algebras comes by precomposing by the action on the operad.  In particular, any framed $E_n$-algebra\footnote{Note the potentially confusing yet standard terminology here: a framed $E_n$-algebra is one which gives invariants of oriented manifolds,} is a homotopy fixed point for the $\mathrm{SO}(n)$ action in a canonical way.  This action is quite special, since it only involves algebra maps and not Morita equivalences.  Similarly the canonical homotopy fixed point structures again only involve algebra maps and not Morita equivalences.  In particular, not all the homotopy fixed points in the Morita category are the canonical ones attached to a framed $E_n$-algebra.  Nonetheless, those special fixed points give oriented TFTs which are sufficient in applications.  

In the $E_2$ setting, a framed $E_2$-algebra is the same thing as a balanced braided tensor category, and so any balanced braided tensor category has a canonical $\mathrm{SO}(2)$-fixed point structure.  This is enough to get a $2$-dimensional oriented TFT, as in~\cite{Ben-Zvi2015}.  But for the $3$-dimensional TFT, having an $\mathrm{SO}(2)$ fixed point structure only yields a ``combed TFT", i.e. one that depends on a choice of non-vanishing vector field on the 3-manifold.

It is relatively easy to characterize $\mathrm{SO}(2)$ actions and their fixed points owing to the standard cell decomposition of $BSO(2) = \mathbb{C}P^\infty$, with a unique cell in each even degree, and with the usual attaching maps.  Since in our case $X$ is the core of a discrete $4$-category $\RigidBrTens$, the target is a homotopy 5-type and by standard arguments we may restrict attention to $5$-type of $\mathbb{C}P^6$. The 2-cell corresponds to the Serre automorphism, which is the value of the loop bordism (an interval with framing that twists once).  From the above description of the $\mathrm{SO}(2)$-action we see that Serre automorphism is given by the double braiding bimodule (see \cite{Ben-Zvi2016}, Figure~3).  The 4- and 6-cells correspond to certain higher compatibilities that the Serre satisfies, but we will not need to understand them in detail for what follows.  Trivializing this action corresponds to picking a trivialization of the Serre, plus some data attached to the higher cells.  A balancing is exactly a trivialization of the Serre that comes from a homomorphism instead of an arbitrary bimodule.  The discussion in the previous paragraph guarantees the existence of a canonical choice for the higher pieces of data.

Since $\pi_1(SO(3)) = \mathbb{Z}/2\mathbb{Z}$, we see that $B\mathrm{SO}(3)$ has an extra $3$-cell which trivializes the square of the Serre.  This is called the belt bordism in~\cite{Douglas2013} because it comes from the belt trick, and its image under the TFT is called the Radford because of its relationship to Radford's theorem for finite tensor categories \cite{Radford1976,Etingof2004a}.  In the case of $\BrTens$ the image of the belt bordism is given by the Drinfeld trivialization of the quadruple braiding bimodule \cite{Drinfeld1989b} \cite[Chapter XIV]{Kassel1995} \cite[\S2.2]{Bakalov2001} \cite[\S8.9.]{Etingof2015}.  So, in order to give a trivialization of the $SO(3)$-action, we must identify the trivialization of the quadruple braiding bimodule given by the Drinfeld map with the one given by the square of the balancing.  In general such an identification is given by a natural transformation of bimodule functors.  But since these bimodules come from tensor functors, and the bimodule functors come from tensor natural isomorphisms, a particularly nice way to give a bimodule natural isomorphism is to simply assert that the two monoidal natural transformations are equal.  This recovers exactly the definition of a ribbon tensor category.  So we see that a ribbon tensor category has a trivialization of the $\mathrm{SO}(2)$-action and a trivialization of the belt bordism.  As before, not all homotopy fixed points will be of this form, but these special homotopy fixed points are sufficient for applications.  (See \cite[\S3.5]{Douglas2013} for analogous partial results.)

We expect it to follow from the main results of \cite{Douglas} that in fact a ribbon category has a canonical $SO(3)$-fixed point structure.  This requires understanding the $5$-type of $BSO(3)$ and checking that certain higher cells automatically vanish because we've chosen the trivialization of the Serre, the higher $SO(2)$-cells, and Radford to be in a very special form (e.g. coming from algebra maps rather than from bimodules).  In a sense this is easier than the corresponding statement about finite tensor categories, because the higher $SO(2)$-cells are easier to trivialize here.

We also speculate that a ribbon fusion category has a canonical $SO(4)$-fixed point structure, but this result would not follow from the results of \cite{Douglas}, but rather would involve computations with the 5-type of $BSO(4)$.  Nonetheless this result should be true by comparison with \cite{Walker} (at least in the unitary case).

\subsubsection{The quantum character field theory}
Our primary motivating application for Theorem~\ref{thm:BrTens-Informal} is to producing a local 3-dimensional topological field theory extending the local 2-dimensional \emph{quantum character theory}, which was introduced in~\cite{Ben-Zvi2015}.  Let us outline some expected applications of this extension here.

Recall that the \emph{character stack} $\uch_G(X)$ associated to a topological space $X$ and a group $G$ is a moduli stack of $G$-local systems on $X$, equivalently, of representations
$$\pi_1(X)\to G,$$
considered up to the conjugation action of $G$.  It is proved in \cite{Ben-Zvi2010} that the assignment to $X$ of its $G$-character stack, accessed through its category of quasi-coherent sheaves, defines an $n$-dimensional topological field theory for any $n$, and that this coincides with the $n$-dimensional field theory assigning $\Rep(G)$ to a point, where we may regard the $E_\infty$-algebra $\Rep(G)$ as an $E_n$ algebra for any $n$.  These may then be computed using the machinery of factorization homology~\cite{Ayala2015,Ayala2017,Beilinson2004,Ginot2015,Lurie}, which has blossomed in recent years.

In low dimensions, character stack exhibit rich additional structures.  It is a celebrated result of Atiyah--Bott~\cite{Atiyah1983} and Goldman~\cite{Goldman1984} that when the group $G$ is reductive, the underlying variety of the $G$-character stack of closed surfaces carry canonical symplectic forms, built using the Killing form for $G$ and the Poincare pairing for $S$.  Furthermore, any 3-manifold $M$ with boundary defines a Lagrangian subvariety $\ch_G(M) \subset \ch_G(\partial M)$ with respect to this symplectic form. This result has been upgraded in the influential papers~\cite{Calaque2013,Calaque2017,Pantev2013} to a so-called 0-shifted symplectic structure on the stack $\uch_G(S)$ for any closed surface $S$, and to a Lagrangian structure on the map
\[
	\uch_G(M)\rightarrow \uch_G(\partial M).
	\]

In~\cite{Ben-Zvi2015}, \emph{deformation quantizations} of the classical character stack were introduced, more precisely deformation quantizations of their categories of quasi-coherent sheaves.  The construction hinged on the machinery of factorization homology, and as such the assignment $S\mapsto QCh_G(S)$ defined a local 2-dimensional topological field theory.  These invariants were computed explicitly for arbitrary punctured surfaces, where they were related to well-known ad-hoc quantizations of character varieties, most notably the moduli algebras of Alekseev~\cite{Alekseev1993,Alekseev1996}.

In the follow-up paper~\cite{Ben-Zvi2016}, the resulting invariants of closed surfaces were computed. Of special interest is the case of the closed torus $T^2$.  The quantum $G$-character stack of $T^2$ was identified with the category of strongly equivariant $\cD_q(G)$-modules.  Here the equivariance condition is with respect to (a quantization of) the conjugation action of $G$ on itself.  This identification allows for the construction of \emph{quantum Hamiltonian reduction} functors, from the quantum character variety of the 2-torus, to an algebra of Weyl group-invariant difference operators $\cD_q(H)^W$ on the Cartan subgroup $H$ of $G$.

The results of the present paper, together with the general formalism of topological field theory, yield several interesting new structures on quantum character varieties, most notably in the case of 3-manifolds with boundary.  Suppose that $M$ is a 3-manifold with a two-torus boundary.  For example, $M$ could be the complement in $S^3$ to a knot $K$.  Then we may regard $M$ as a cobordism from the empty manifold $\emptyset$ to $T^2$, giving rise to a functor,
$$\Vect = QCV(\emptyset) \xrightarrow{QCV(M)} QCV(T^2) \xrightarrow{\textrm{Q.H.R.}} \cD_q(H)^W\modu.$$
Identifying the functor obtained in this way with its value on the one-dimensional vector space $\K\in\Vect$, we obtain from any such $M$ a module for the ring $\cD_q(H)^W$.  The isomorphism class of such a module is a homeomorphism invariant of $M$, which may be regarded as a quantization of the defining Lagrangian of the classical character variety of $M$.  As such, it should provide a foundational framework for:  quantum $A$-polynomials~\cite{Dimofte2013,Garoufalidis2004,Gukov2005,Gukov2012}, DAHA-Jones polynomials and skein modules of knots complements~\cite{Cherednik2013,Berest2016,Samuelson2017}, and refined Chern-Simons invariants~\cite{Aganagic2011}, among others, each of which yet lacks a topologically invariant formulation. 

\subsubsection{Skein modules, skein algebras, and skein categories}\label{sec:skeins}
The Poisson algebra of functions on the classical $G$-character variety on a punctured surface $S$, has a natural description in terms of graphs drawn on the surface, with edges labelled by $G$-modules, and vertices labelled by $G$-morphisms. This led to a successful program of deformation quantization of character varieties via ``skein theory" approach~\cite{Frohman2000,Przytyckia2000,Andersen1998,Roche2002}, which we now recall.

Given an oriented 3-manifold $M$, its \emph{skein module} A(M) is the quotient of the vector space freely generated by isotopy classes of links in $M$, modulo ``skein relations": these relate vectors corresponding to tangles obtained by certain basic modifications within a fixed ball (see e.g.~\cite{Lickorish2009} for a survey).  Perhaps the best-studied of these is the Kauffmann bracket skein relations~\cite{Przytycki2006,Kauffman2001} which derive from the representation theory of $U_q(\mathfrak{sl}_2)$, but there are also the HOMFLYPT polynomials, which derive from the $N\to \infty$ limit of the representation theories of $U_q(\mathfrak{sl}_N)$ and the Kauffman polynomial which similarly corresponds to quantum orthogonal and symplectic families.  More generally, one associates to any 3-manifold and any (small) ribbon tensor category the associated skein module, defined using ribbon graphs labelled by objects and morphisms of $\cA$.

In the case where $M=S\times I$ for some surface $S$, the skein module acquires a natural algebra structure -- the ``skein algebra" -- by concatenating in the interval direction.  In the case of finite dimensional modules over the quantum group of $G$, this algebra quantizes the classical character variety~\cite{Roche2002}. Skein modules for 3-manifolds with boundary naturally inherit the structure of a module over the skein algebra of their boundary in a similar manner.  Allowing surfaces with marked points, one is naturally led to the so-called skein categories $C(S)$ of Morrison-Walker \cite{Walker,Morrison2012,Morrison2011}.  An object of the skein category is a finite collection of points of the surface, colored by objects of the ribbon braided tensor category; Homs are given by \emph{relative skein modules} in the cylinder, i.e. ribbon graphs drawn in $S\times I$ ending at the top and bottom of the cylinder compatibly with their markings.  In the same way, one attaches to an arbitrary oriented 3-manifold $M$ with oriented boundary $\partial M = S_+ \cup S_-$ a functor $C(S_{+})^{op}\times C(S_{-})\to \Vect$ (sometimes called a ``bimodule"), which outputs the skein module in the 3-manifold, relative to the markings on its boundary which define the input objects.

Our constructions in the current paper extend the 2D theory of \cite{Ben-Zvi2015} to a 3D theory; however the factorization homology construction is specific to dimension 2 (or to dimension $n$, more generally, when working with $E_n$-algebras). Thus it is an important and interesting question to describe the resulting functors explicitly, in the framework of skein modules of 3-manifolds.  In the semi-simple case, we outline below a conjectural answer to this question, following a proposal of \cite{Johnson-Freyd2015}.

Recall that in \cite{Ben-Zvi2015} the role of skein categories is played instead by factorization homology with coefficients in a braided tensor category, as the basic ingredient in the construction.  Like skein categories, factorization homology categories feature special objects indexed by points of the surface colored by representations of the input braided tensor category.  Whereas in skein theory, these are taken as a starting point to define the categories, in factorization homology they are instead induced via functoriality by disk embeddings into the surface.

Factorization homology categories may be understood therefore as a suitable co-completion of the skein categories, in that they contain skein categories as a full, generating, subcategory, and they further contain additional objects obtained from those by taking cokernels and direct sums. This relationship is analogous to the obtainment of the category of coherent sheaves on a variety as a co-completion of the category of vector bundles.  

Assuming the orientability conjectures from Section \ref{sec:oriented}, an oriented 3-manifold $M$ with boundary $S^+\sqcup S^-$, viewed as a cobordism, gives rise to a functor $F_M$ between the factorization homology categories.  Because the factorization homology categories are \emph{generated} by the skein-theoretic subcategory, in order to give a functor between factorization homology categories, it suffices to specify its values on this subcategory.
We further conjecture, following \cite{Walker,Johnson-Freyd2015}, that $F_M$ is given on skein objects $X$ and $Y$ by the formula:
$$\Hom(F_M(X),Y) = Sk^{rel}_{X,Y}(M)$$
Here the RHS denotes the relative skein module, spanned by ribbon graphs in $M$ ending at the points compatibly with the labeling of $X$ and $Y$.  Interpreting presentable $\K$-linear categories as categorified Hilbert spaces, the $\Hom$ functor plays the role of categorified inner product, and so this gives a sort of ``matrix coefficient" for the functor.  We note that such a formula can only be expected to hold in the semi-simple case, as it implies that the unit is compact projective.

\subsubsection{Crane--Yetter and Witten-Reshetikhin--Turaev theory}
As a consequence of Theorem~\ref{thm:BrFus-Informal}, every braided fusion category gives rises to a local framed 4d TFT. As explained above, starting from a \emph{ribbon} braided fusion category, also known as a pre-modular category, one conjecturally obtains an oriented theory, hence numerical invariants of oriented 4-manifolds. We expect this invariant to recover the one defined by Crane--Yetter~\cite{Crane1993} in the modular case and by Crane--Kauffman--Yetter~\cite{Crane1997} in the general pre-modular case. 

In the modular case, unpublished work of Freed--Teleman~(see \cite{Freed2012a,Freeda}) and Walker (see \cite{Walker}) asserts that the corresponding local theory is \emph{invertible} and encodes the anomaly of Witten--Reshetikhin--Turaev theory \cite{Reshetikhin1991,Witten1989}. Roughly speaking, this means that the invertible theory can be upgraded to a \emph{relative} field theory~\cite{Freed2012,Fiorenza2015,Johnson-Freyd2017}, which in particular can be used to define a numerical invariant of pairs $(M,[W])$ where $W$ is a 4-manifold with boundary $M$ and $[W]$ is the cobordism class of $W$. This invariant is, in turn, expected to coincide with the Witten--Reshetikhin--Turaev invariant attached to the same data. We note that the fact that Witten--Reshetikhin--Turaev theory can be thought of as living on the boundary of an almost trivial 4-dimensional theory is apparent from its relation with Chern--Simons theory~\cite{Witten1989}, and is essentially the approach outlined in~\cite{Walker}.  

Moving beyond the modular case, non-semisimple versions of the Witten-Reshetikhin-Turaev invariants have recently been constructed in a series of papers~\cite{Beliakova,Blanchet,DeRenzi}.  It is natural to expect that the TFT determined by a finite and factorizable, but non-semisimple braided tensor category, regarded as an object of $\BrTens$, can be related to those in a similar way.  We hope to return to this in future work.

\subsection{Outline}
Section~\ref{sec:prelim} contains a recollection on locally presentable tensor and braided tensor categories, as well as on the notion of dualizability in higher categories. In Section~\ref{sec:BrTens} we give a detailed description of the symmetric monoidal Morita $4$-category of braided tensor categories. Section~\ref{sec:rigidity} is devoted to the definition of rigidity in our framework, and the proof that restricting to cp-rigid tensor and braided tensor categories, and bimodules with enough compact-projective over those, form a sub-category of the Morita category. Finally, Section~\ref{sec:dualizability} contains the proof of some dualizability statements which are implied by cp-rigidity, and finally the proof our main results.
\subsection{Acknowledgements}
AB is supported by the RTG 1670 ``Mathematics inspired by String Theory and Quantum Field Theory''. DJ is supported by European Research
Council (ERC) under the European Union's Horizon 2020 research and innovation
programme (grant agreement no. 637618).
NS is supported by NSF grant DMS-1454767. 

We would like to thank David Ben-Zvi, Christopher Douglas, Theo Johnson-Freyd, Tom Leinster,  Claudia Scheimbauer, Christopher Schommer-Pries, Mike Shulman, and Kevin Walker for helpful conversations.  We would particularly like to thank Theo Johnson-Freyd for some very informative discussions about different definitions of rigidity, and how they resolve some pathologies in characteristic $p$, and David Ben-Zvi for sharing his vision of Betti geometric Langlands and its quantization, which motivates this work, as well as helpful discussions about rigidity in the $\oo$-category setting.

\section{Preliminaries}\label{sec:prelim}

In this section, we recall a number of basic definitions we will need, and we give a more complete statement of our main results.  Let us fix throughout a field $\K$ of characteristic zero. General references include~\cite{Adamek1994,Brandenburg2015,Makkai1989}.

\subsection{Presentable linear categories}

\begin{definition}
A category is said to be $\K$-linear if it is enriched and tensored over the category of $\K$-linear vector spaces.\end{definition}

\begin{definition}
Let $\cC$ be a $\K$-linear category. An object $c \in \cC$ is called
\begin{itemize}
	\item \emph{compact}, if $\Hom(c,-)$ commutes with filtered colimits
	\item \emph{compact-projective}, if $\Hom(c,-)$ is cocontinuous, i.e. commutes with arbitrary small colimits.
		\end{itemize}
	\end{definition}

\begin{definition}
We say that a $\K$-linear category
\begin{itemize}
\item is \emph{locally finitely presentable} if it admits arbitrary small colimits, and is generated under filtered colimits by a small sub-category of compact objects. 
\item \emph{has enough compact projectives} if it admits arbitrary small colimits, and is generated under small colimits by a small sub-category of compact-projective objects.
\end{itemize}
\end{definition}

\begin{remark} We note that, by a standard exercise, a category which has enough compact projectives in the above sense is in particular locally finitely presentable.\end{remark}

\begin{remark}\label{rem:projective}
In the $\K$-linear setting an object is compact-projective if, and only if, it is compact and is projective in the sense that $\Hom(p,-)$ is right exact (preserves finite colimits).  As usual, projective has many equivalent formulations, including one that says $p$ is projective if for any morphism $f: p\rightarrow b$ and any epimorphism $e: a \rightarrow b$ then $f$ factors through $e$.  We will need this equivalent formulation in Lemma \ref{lem:absoluteprojectives}.  Note that in general saying every SES $0 \rightarrow a \rightarrow b \rightarrow p \rightarrow 0$ splits is necessary but not sufficient for $p$ to be projetive, though it is sufficient for abelian categories.  See Qiaochu Yuan's blog posts for clear exposition of these standard results in the cocomplete setting \cite{YuanTiny,YuanProjective}.
\end{remark}

\begin{remark}
The phrase `$\cC$ has enough projectives' is typically used to mean that every object $X$ of $\cC$ has a projective cover.  In that case, choosing first a projective cover $p: P\to X$, and then a projective cover $p':P'\to \ker(p)$, we can write $X\cong\coker(p')$ as colimit of projectives.  This motivates the terminology `enough compact projectives'.  Compact objects are also called ``small" and compact projectives are also called ``small projective" or``tiny" \cite{Kelly1982, Yetter1987, Brandenburg2015}.
\end{remark}
	
In enriched category theory, it is standard to consider $\lambda$-locally presentable categories for a regular cardinal $\lambda$ (see \cite{Adamek1994} for details); taken together these form the \emph{locally presentable} categories.  Locally finite presentable categories correspond to $\lambda=\aleph_0$, while categories with enough compact projectives in the above sense correspond to $\lambda=0$ \cite{Brandenburg2015}.  For a general regular cardinal $\lambda$, $\lambda$-locally presentable categories are defined similarly to the locally finitely presentable categories, except that the cardinal $\lambda$ enters into the notion of compact objects (their Hom functors preserve only $\lambda$-filtered objects), and the notion of generation (the $\lambda$-compact objects generate under $\lambda$-filtered colimits).  Since we will not require any technical arguments involving general cardinals, the reader is safe to think of locally finitely presentable categories in all arguments, without loss of generality.  

\begin{definition} We denote by:
\begin{itemize}
	\item $\LFP$ the 2-category of locally presentable categories, cocontinuous functors and natural transformations\footnote{We note that this terminology differs slightly from \cite{Ben-Zvi2015}, where only natural isomorphisms were allowed},
\item $\LFPd$ the full 2-subcategory consisting of categories having enough compact projectives.
\end{itemize}
\end{definition}

The 2-category of $\K$-linear categories is symmetric monoidal:  given two $\K$-linear categories, their linear tensor product $\cC\ot\cD$ has as its objects pairs of objects of $\cC$ and $\cD$, and morphisms defined by
\[
	\Hom_{\cC\ot \cD}((c_1,d_1),(c_2,d_2)):=\Hom_\cC(c_1,c_2)\ot_\K \Hom_\cD(d_1,d_2).
	\]
However, the $\K$-linear tensor product of two categories in $\LFP$ is not again in $\LFP$.  Each of $\LFP$ and $\LFPd$ nevertheless admit a natural symmetric monoidal structure, extending the $\K$-linear tensor product, and define as follows:

\begin{definition}
	The Deligne--Kelly~\cite{Kelly1982, Deligne2007, Franco2013} tensor product of categories $\cC,\cD \in \LFP$ is another category $\cC\bt \cD\in\LFP$, equipped with a linear functor $\cC \otimes \cD \rightarrow \cC\bt \cD$, cocontinuous in each variable, which is moreover universal for this property, in the sense that we have an equivalence of groupoids: 
	$$\Hom_{\LFP}(\cC\boxtimes \cD,\cE)^\times \simeq \Lin^{c,c}(\cC\otimes\cD,\cE)$$
\end{definition}
\noindent where on the left hand side we throw away non-invertible natural transformations and $\Lin^{c,c}$ denotes the groupoid whose objects are $\K$-linear functors which are co-continuous in each variable, and whose morphisms are natural isomorphisms.

\begin{proposition} \label{prop:DKgen}
The Deligne--Kelly tensor product $\cC\bt \cD$ exists and is again locally presentable.
\end{proposition}
\begin{proof} It is shown in \cite[Lemma 8]{Franco2013}, following \cite[Chapter 6]{Kelly1982} (see also \cite{Caviglia}) that the Deligne--Kelly tensor product of $\lambda$-presentable categories $\cC$ and $\cD$ is again $\lambda$-presentable, and is in fact generated under $\lambda$-filtered colimits by $\lambda$-small colimits of pure tensor products of the $\lambda$-compact generators in each component. \end{proof}
\begin{corollary}\label{prop:purecp}
The Deligne-Kelly tensor product of two categories with enough compact projectives again has enough compact projectives, and is generated by the pure tensor products of the compact projective generators in each component.
\end{corollary}

\begin{remark} \label{rem:Rex}
	Most examples considered in the tensor category literature (e.g. \cite{Etingof2015}) are not in $\LFP$, because infinite direct sums are excluded.  For example, such literature more commonly considers finite dimensional vector spaces rather than all vector spaces.  However, these examples all fit into the setting of finitely cocomplete categories (called $\Rex$), and since the Ind-completion of a $\Rex$ category is $\LFP$ while the subcategory of compact objects in an $\LFP$ category is $\Rex$ one can easily translate any of these examples into our setting.  In particular, the Ind-completion of any locally finite category (abelian, finite dimensional Hom-spaces, finite length, and enough projectives) is $\LFP$.  However, if one works directly with finitely cocomplete categories and their right exact functors this is insufficient to get our dualizability results: under Ind-completion right exact functors corresponds to \emph{compact preserving} cocontinuous functors, a property which will typically fail for coevaluation maps. Instead one needs to work with finitely cocomplete categories and their right exact \emph{profunctors} (i.e. functors to the Ind-completions).  The importance of using profunctors (or bimodules) instead of functors in order to get stronger dualizability results also appears in \cite{Stay2016} and \cite{Walker}.  See \cite[\S3]{Ben-Zvi2015} for more details on the relationship between $\Rex$ and $\LFP$.  
\end{remark}

\subsection{Dualizability}
We recall the following from \cite{Lurie}. Roughly speaking, a symmetric monoidal $n$-category is said to be $k$-dualizable for some $1\leq k\leq n$ if every object has a dual in the usual sense, and every $i$-morphism for $1\leq i < k$ has a left and a right adjoint. This in particular implies 
that every object is $k$-dualizable in the sense that its evaluation and coevaluation maps have left and right adjoint, the units and counits of which themselves have adjoints, and so on.

This is made more precise as follows. Let $\cC$ be an $n$-category and let $h_2\cC$ be its homotopy 2-category, i.e. the 2-category whose objects and morphisms are the same as $\cC$ and whose 2-morphisms are isomorphism classes of 2-morphisms of $\cC$.
\begin{definition}
Let $\cC$ be a 2-category and $x,y\in \cC$. Two morphisms $L:x\rightarrow y$ and $R:y\rightarrow x$ are said to be adjoints if there exists 2-morphisms $\eta:\id_x\rightarrow R\circ L$ and $\epsilon :L\circ R\rightarrow \id_y$ such that
	\begin{align*}
		(R\epsilon)\cdot (\eta R)&=\id_R&(\epsilon L)\cdot (L \eta)&=\id_L
	\end{align*}
\end{definition}

\begin{example}
The motivating example for the above definition is of course the $2$-category of categories, functors and their natural transformations, where this recovers the usual notion of left and right adjoints of functors.  Another important example related to rigidity comes by regarding a monoidal category $\cC$ as a 2-category $B\cC$, with a unique object $\ast$, with $End(\ast):=\cC$.  Then an object of $\cC$, regarded as a 1-morphism of $B\cC$, has an left/right adjoint if and only if that object has a left/right dual as an object of $\cC$.
\end{example}

\begin{definition}  
	Let $\cC$ be an $n$-category and let $1\leq k < n$. The category $\cC$ is said to have adjoints for $k$-morphisms if it satisfies the following inductively defined property:
	\begin{itemize}
		\item if $k=1$, then in $h_2 \cC$ every morphism has a left and a right adjoint
		\item if $k>1$, then for every pair of objects $x,y \in \cC$, the $(n-1)$-category $\Hom(x,y)$ has adjoints for $k-1$-morphism.
	\end{itemize}
\end{definition}
Let now $\cS$ be a symmetric monoidal $n$-category and define $B\cS$ to be the $(n+1)$-category with one object $*$ and morphisms given by $\End(*):=\cS$.
\begin{definition}
We say the category $\cS$ is $k$-dualizable if the category $B\cS$ has adjoints for $k$-morphisms.
\end{definition}
\begin{remark}
We stress that the notion of adjoints in a (possibly non-discrete) $2$-category requires the zig-zag identities to hold only up to \emph{non-specified} higher isomorphisms. In other words this is a property, not extra structure, which is why it is enough to check this in the (discrete) homotopy 2-category. In particular, although the notion of being $k$-dualizable is typically higher categorical, unravelling the definition it boils down to check a list of conditions in various discrete 2-categories.  See the remark in Section I.5 of \cite{Douglas2013}.
\end{remark}

\begin{remark}\label{leftright}
We will use the same left/right conventions for duality and adjoints as in \cite[\S2.1]{Douglas2013}.  In particular, a bimodule ${}_{\cC}\cM_{\cD}$ is regarded as a morphism from $\cC$ to $\cD$, and so a left adjoint of a bimodule ${}_{\cC}\cM_{\cD}$ is a bimodule ${}_{\cD} \cN_{\cC}$ with an evaluation map $\ev:  {}_{\cD} \cN \boxtimes_{\cC} \cM_{\cC} \rightarrow {}_\cD \cD_\cD$ and a coevaluation $\coev:  {}_\cC \cC_\cC \rightarrow {}_{\cC} \cM \boxtimes_{\cD} \cM_{\cC}$ satisfying the usual zig-zag.
\end{remark}

\subsection{Tensor categories}
The following definitions are standard, but we emphasize that throughout this paper we always work in $\LFP$.  For any details about what ``the obvious diagrams'' means in some of these definitions see \cite{Etingof2015}.

\begin{definition} A \emph{tensor category} is a category $\cA$ in $\LFP$ equipped with a cocontinuous functor (i.e. a morphism in $\LFP$)
	\[
\ot: \cA \bt \cA \longrightarrow \cA,
	\]
	a distinguished object $\un_\cA$ (the unit) and natural isomorphism
	\[
		\alpha_{x,y,z}:(x\ot y)\ot z \longrightarrow x \ot (y\ot z)
	\]
	and
	\begin{align*}
		l_x:\un_\cA\ot x &\rightarrow x & r_x:x \ot \un_\cA \rightarrow x
	\end{align*}
	such that the following diagram commutes:
\begin{center}\begin{tikzpicture}[description/.style={fill=white,inner sep=2pt}]
\matrix (m) [matrix of math nodes, row sep=2em,
column sep=-0.5em, text height=1.5ex, text depth=0.25ex]
{&&(x\otimes(y\otimes z))\otimes w&& \\
((x\otimes y)\otimes z)\otimes w&&&&x\otimes((y\otimes z)\otimes w)\\
    (x\otimes y)\otimes(z\otimes w)&&&&x\otimes(y\otimes (z\otimes w))\\};
\path[->,font=\scriptsize]
(m-2-1.north) edge node[auto] {$ \alpha_{x,y,z} \ot \id_w $} (m-1-3.south west)
(m-1-3.south east) edge node[auto] {$\alpha_{x,y\ot z,w}$} (m-2-5.north)
(m-2-5.south) edge node[auto] {$ \id_x \ot \alpha_{y,z,w} $} (m-3-5.north)
(m-3-1.east) edge node[below] {$ \alpha_{x,y,z\ot w} $} (m-3-5.west)
(m-2-1.south) edge node[left] {$ \alpha_{x\ot y,z,w} $} (m-3-1.north)
;
\end{tikzpicture}\end{center}

\begin{center}\begin{tikzpicture}[description/.style={fill=white,inner sep=2pt}]
\matrix (m) [matrix of math nodes, row sep=2em,
column sep=0.5em, text height=1.5ex, text depth=0.25ex]
{        (x \ot \un_\cA) \ot y && x \ot (\un_\cA \ot y)\\
			       &x \ot y\\
};
\path[->,font=\scriptsize]
(m-1-1) edge node[auto] {$ \alpha_{x,\un_\cA,y} $} (m-1-3)
(m-1-1) edge node[below] {$ r_x $} (m-2-2)
(m-1-3) edge node[auto] {$ l_x $} (m-2-2)
;
\end{tikzpicture}.\end{center}

\end{definition}
\begin{remark}Equivalently, a tensor category is an $E_1$-algebra in $\LFP$.
\end{remark}
\begin{remark}
As customary, we suppress $\alpha$, $l$ and $r$ when those are clear from the context, as they can be uniquely filled in.
\end{remark}
\begin{definition}
	Given a tensor category $\cC$, we define the \emph{multiplication opposite} $\cC^\mop$ as the tensor category with reversed multiplication and inverse associativity constraint.
\end{definition}
\begin{definition}
Let $\cA$ be a tensor category. A left (resp, right) $\cA$-module is a category $\cM\in \LFP$ equipped with a cocontinuous functor, $\ot:\cA \bt \cM\rightarrow \cM$ (resp, $\ot:\cM \bt \cA \rightarrow \cM$), together with an associativity constraint and a natural isomorphism $\un_\cA\ot m\rightarrow m$ for $m\in \cM$ making the analogous pentagon diagram commute.
\end{definition}
\begin{definition}
Let $\cA,\cB$ be tensor categories. An $\cA-\cB$-bimodule is an $\LFP$ category $\cM$ which is simultaneously a left $\cA$-module and a right $\cB$-module, together with an isomorphism
\[
	\Gamma:(a\ot m)\ot b \longrightarrow a \ot (m \ot b)
\]
making the obvious diagram commute.
\end{definition}

\begin{definition}
	Let $\cM$, $\cN$ be two module categories over $\cA$. Then a \emph{left module functor}, or $\cA$-linear functor, is a pair of a functor $F:\cM\rightarrow \cN$ and a natural isomorphism $f:F(a\ot m)\rightarrow a\ot F(m)$ for $a \in \cA,m\in \cM$ making the obvious diagrams commute. Right module and bimodule functors are defined similarly.
\end{definition}
\begin{definition}
Let $\cA$ be a tensor category and $\cM,\cN$ be a right and a left $\cA$-module category respectively. An $\cA$-balanced functor is a pair of a functor from $F:\cM \bt \cN$ to some $\LFP$ category $\cE$ and a natural transformation 
	\[
		f:F((m\ot a)\bt n)\cong F(m\bt (a \ot n))
		\]
		for $a\in \cA,m\in\cM,n\in \cN$ making the obvious diagrams commute. 
\end{definition}
\begin{definition}
Let $\cA$ be a tensor category and let $\cM$ and $\cN$ be a right and left $\cA$-module, respectively. The balanced (or relative) Deligne--Kelly tensor product is another category $\cM\bt_\cA \cN\in\LFP$ together with an $\cA$-balanced functor $\bt_\cA:\cM\bt \cN \rightarrow \cM \bt_\cA \cN$ which induces a natural equivalence of categories between balanced functor out of $\cM \bt \cN$ to some category $\cE$, and morphisms in $\LFP$ from $\cM \bt_\cA \cN$ to $\cE$.  The existence of the balanced Deligne-Kelly tensor product follows from the cocompleteness of $\LFP$, see \cite[Def. 3.13, Rem. 3.14]{Ben-Zvi2015}.  Constructions of the balanced tensor product in special cases appear in \cite{Etingof2010, Davydov2013, Douglas2014}.
\end{definition}

\begin{definition}
A \emph{fusion category} is a cp-rigid $\LFP$-tensor category, which is semi-simple (in the sense that every object is a possibly infinite direct sum of simple objects), with simple unit, and with finitely many isomorphism classes of simple objects.
		\end{definition}
		\begin{remark}
			More precisely, what we define here is the ind-completion of what is called a fusion category in the literature.  See Remark \ref{rem:Rex}.  Note that under the semi-simplicity assumption every compact object is projective, so the notions of compact-rigid and cp-rigid coincide.
		\end{remark}
\subsection{Braided tensor categories}

\begin{definition}\label{defi:BTC}
	A braided tensor category is a tensor category $(\cA,\ot,\alpha)$ together with a natural automorphism $\beta$ of $- \ot -$ making the following diagrams commute:
\begin{center}\begin{tikzpicture}[description/.style={fill=white,inner sep=2pt}]
\matrix (m) [matrix of math nodes, row sep=3em,
column sep=2.5em, text height=1.5ex, text depth=0.25ex]
{& x \ot (y \ot z) & (y \ot z) \ot x &\\
(x \ot y) \ot z & && y \ot (z \ot x) \\
&(y \ot x) \ot z & y \ot (x \ot z)& \\};
\path[->,font=\scriptsize]
(m-1-2) edge node[auto] {$\beta_{x,y\ot z}$}(m-1-3)
(m-1-3) edge node[auto] {$\alpha_{y,z,x}$}(m-2-4)
(m-3-3) edge node[below right] {$\id_y\ot \beta_{x,z}$}(m-2-4)
(m-3-2) edge node[below] {$\alpha_{y,x,z}$}(m-3-3)
(m-2-1) edge node[below left] {$\beta_{x,y}\ot \id_z$}(m-3-2)
(m-2-1) edge node[auto] {$\alpha_{x,y,z}$}(m-1-2)
;
\end{tikzpicture}\end{center}
\begin{center}\begin{tikzpicture}[description/.style={fill=white,inner sep=2pt}]
\matrix (m) [matrix of math nodes, row sep=3em,
column sep=2.5em, text height=1.5ex, text depth=0.25ex]
{ & (x \ot y) \ot z  & z \ot (x \ot y) \\
x \ot (y \ot z)  & && (z \ot x )\ot y \\
&x \ot (z \ot y) & (x \ot z) \ot y 
\\};
\path[->,font=\scriptsize]
(m-1-2) edge node[auto] {$\beta_{x\ot y,z}$}(m-1-3)
(m-1-3) edge node[auto] {$\alpha_{z,x,y}^{-1}$}(m-2-4)
(m-3-3) edge node[below right] {$ \beta_{x,z}\ot \id_y$}(m-2-4)
(m-3-2) edge node[below] {$\alpha_{x,z,y}^{-1}$}(m-3-3)
(m-2-1) edge node[below left] {$\id_x \ot \beta_{y,z}$}(m-3-2)
(m-2-1) edge node[auto] {$\alpha_{x,y,z}^{-1}$}(m-1-2)
;
\end{tikzpicture}\end{center}

\end{definition}
\begin{remark}
Equivalently, a braided tensor category is an $E_2$-algebra in $\LFP$.
\end{remark}

\begin{definition} Let $\cA$ be a braided tensor category.  We define \emph{braiding reverse} of $\cA$ to be the braided tensor category $\cA^{\bop}$ to having the same underlying tensor category, but with braiding isomorphism $\sigma_{V,W}$ replaced by $\sigma_{W,V}^{-1}$.
\end{definition}

\begin{remark}
The opposite braiding comes from reflecting the discs about the $x$-axis in the $E_2$ operad and can be thought of as the opposite in the second multiplication direction.  The opposite in the first multiplication direction is $\cA^\mop$ with the reversed tensor product and the braiding given by $\sigma_{V,W}^{-1}: V \otimes^{op} W = W \otimes V \rightarrow V \otimes W = W \otimes^{op} V$.  Again $\cA^\mop$ corresponds to a reflection, this time about the $y$-axis.  It is not difficult to see that $\cA^\mop$ and $\cA^\bop$ are isomorphic using the braiding.  Note that the double opposite $\cA^{\mop\bop} = \cA^{\bop\mop}$ has underlying tensor category $\cA^\mop$ but with the braiding given by $\sigma_{W,V}$.  Since $\cA^{\mop\bop}$ corresponds to rotation by $180$-degrees it is orientation preserving and should not be thought of as an opposite (it is isomorphic to the original $\cA$ and not to either of the opposites).  
\end{remark}

Finally, let us recall the following well-known construction of a braided tensor category out of a tensor category.
\begin{definition}
	Let $(\cC,\ot,\alpha)$ be a tensor category. Then its \emph{Drinfeld center}, or simply \emph{center}, is a braided tensor category $Z(\cC)$ defined as follows
	\begin{itemize}
\item objects are pairs $(y,\beta)$ where $\beta$ is a natural isomorphism
	\[
\beta_x:x \ot y \longrightarrow y \ot x
	\]
	making the obvious analog of the second diagram in Definition~\ref{defi:BTC} commutes.
\item a morphism $(y,\beta) \rightarrow (y',\beta')$ is a morphism $f:y\rightarrow y'$ such that
	\[
		\forall x\in \cC,\ (f \ot \id_x)\beta_x=\beta'_x( \id_{x}\ot f)
	\]
\item 	the tensor product of $(y,\beta)$ and $(y',\beta')$ is the pair ($y\ot y',\tilde \beta)$ where $\tilde\beta$ is defined by the first diagram of Definition~\ref{defi:BTC} with $y'$ instead of $z$
\item the braiding of $(y,\beta)\ot (y',\beta')$ is simply given by $\beta'_y$.
	\end{itemize}
\end{definition}
\begin{remark}
Once again this turns out to be a particular case of the general formalism of $E_n$-algebras: to any $E_n$-algebra in a sufficiently nice symmetric monoidal $1$-category (typically, non-discrete) one associates its  Hochschild cohomology (also called its center), which has a natural structure of an $E_{n+1}$-algebra. 
\end{remark}
The following properties of the center are straightforward and well-known:
\begin{proposition}\label{prop:centerProperty1}
Let $\cC$ be a tensor category. The assignment $(y,\beta)\mapsto (y,\beta^{-1})$ induces a braided tensor equivalence
	\[
		Z(\cC^\mop)\longrightarrow Z(\cC)^{\bop}.
	\]
\end{proposition}
\begin{proposition}\label{prop:centerProperty2}
Let $\cA$ be a braided tensor category with braiding $\beta$. Then there are braided tensor functors
	\begin{align*}
		\cA &\longrightarrow Z(\cA) & \cA^{\bop}& \longrightarrow Z(\cA)\\
		x &\longmapsto (x,\beta_{-,x})& x &\longmapsto (x, \beta_{x,-}^{-1})
	\end{align*}
	which assemble into a single braided tensor functor $\cA \bt \cA^{\bop} \rightarrow Z(\cA)$.
\end{proposition}

It is a general fact that the category of modules over an $E_2$-algebra is an $E_1$, i.e. monoidal, category. Specializing in the case at hand this recovers the following well-known
\begin{proposition}\label{prop:AmodMonoidal}
Let $\cA$ be a braided tensor category with braiding $\beta$. 
	\begin{itemize}
		\item Every left $\cA$-module category $\cM$ (which we assume to be strict) has a canonical structure of a right $\cA$-module, with the same action, and associativity constraint given, for $a,b\in \cA,m\in \cM$
			\[
				(a \ot b)\ot m \xrightarrow{\beta_{a,b}} b \ot (a \ot m).
				\]

\item Given two left $\cA$-modules $\cM,\cN$, the balanced tensor product $\cM \bt_\cA \cN$ where $\cM$ is given the above right module structure, turns the category of left $\cA$-modules into a monoidal 2-category.
	\end{itemize}
\end{proposition}
\section{The Morita theory $\BrTens$}\label{sec:BrTens}

In this section we describe in detail objects and the 1-, 2-, 3- and 4-morphisms which comprise the $4$-category $\BrTens$, as well as their various compositions.  We do not treat in the same detail the simpler case of $\Tens$, partly because this case is already well-documented~\cite{Douglas2013,Etingof2010,Greenough2010,Drinfeld2010}, and partly because a complete definition of $\Tens$ may be recovered from $\BrTens$ by considering the endomorphism $3$-category of the unit object $\Vect\in\BrTens$.

\subsection{Models for higher Morita theories}
There are two closely related models for higher Morita theories of $E_n$-algebras.  Following a proposal of Lurie in~\cite{Lurie2009}, Scheimbauer gave a geometric construction in \cite{Scheimbauer2014}, using the framework of locally constant factorization algebras on stratified spaces of the form $(\RR^k\subset\RR^{k+1}\subset\cdots\subset\RR^n)$ to build an $(n+1)$-category of $E_n$-algebras. Independently, Haugseng~\cite{Haugseng2017} gave an alternative combinatorial/operadic construction of an $(n+1)$-category of $E_n$-algebras, building explicit simplicial sets capturing the basic structures, operations, and coherences capturing the iterative notion of bimodules for $E_n$-algebras.

While there are a number of expected relations between the two constructions, there are not at present writing any theorems relating them, and it appears to be a formidable undertaking to construct such a relation precisely.  At the level of their objects, such a relation is Lurie's theorem asserting the equivalence between $E_n$-algebras and locally constant factorization algebras on $\mathbb{R}^n$ (see \cite[\S 5.2.4]{Lurie2009}, and see \cite{Ginot2015} for an exposition).  However extending Lurie's equivalence to the entire Morita theory seems to be a very difficult task.

An important distinction between the two models for higher Morita theories is the role of \emph{pointings}, i.e. maps from the unit of $\cS$ to all outputs in the theory.  The geometric formulation of Scheimbauer's constructions endows objects and morphisms at all levels with pointings; essentially these arise because the empty open set is initial in the category of open sets on any manifold, so it induces a canonical map to the value of the factorization algebra on any open set, compatibly with push-forward operations.  These pointings are not present in Haugseng's work, except in the low degrees where they endow $E_k$-algebras in the theory with units.

Subsequently to their introduction, each of these two higher Morita theories was extended in \cite{Johnson-Freyd2017} to so-called ``even higher Morita theories", in which the symmetric monoidal $1$-category $\cS$ is replaced by a symmetric monoidal $k$-category.  Correspondingly, instead of an $(n+1)$-category of $E_n$-algebras, one obtains an $(n+k)$-category of $E_n$-algebras.  One of the main examples highlighted in~\cite[Example 8.9]{Johnson-Freyd2017} is precisely the one we need in the present work: the $4$-category $\BrTens$, of braided tensor categories, viewed as $E_2$-algebras in the $2$-category $\LFP$.

Our theorems establishing dualizability in $\Tens$ and $\BrTens$ hold only in the unpointed formulation of Haugseng:  in the pointed setting it is shown in \cite{Gwilliam2018} that $E_n$ algebras are no more than $n$-dualizable.  For this reason, we deploy Haugseng's model in our constructions in this section.  Nonetheless, each structure appearing in the construction of $\BrTens$ can indeed be understood in the language of locally constant factorization algebras, and this gives useful geometric motivation for often elaborately phrased categorical definitions.  For this reason, we also outline in Section~\ref{sec:lcfa} the relationship to Scheimbauer's Morita theory of locally constant factorization algebras, at an informal level and simply ignoring the issue of pointings.

An important feature of Haugseng's construction is Theorem 5.49 from~\cite{Haugseng2017}, which states that the $\Hom$ $n$-categories between two $E_n$-algebra $A$ and $B$ can be identified inductively with the $n$-category of $E_{n-1}$-algebras in the $E_{n-1}$-category of $A$-$B$-bimodules.  Applying this description inductively, in combination with \cite{Johnson-Freyd2017} in the top degrees, we obtain the following description of the Morita theory of $E_2$-algebras:

\begin{itemize}
\item Objects are given by of $E_2$-algebras $\cA$, $\cB$, \ldots
\item $\Hom(\cA,\cB)$ consists of $E_1$-algebras $\cC$,$\cD$, \ldots \emph{internal to} $\cA$-$\cB$-bimodules.
\item $\Hom(\cC,\cD)$ consists of $\cC$-$\cD$-bimodules $\cM$, $\cN$, \ldots \emph{internal to} $\cA$-$\cB$-bimodules.
\item $\Hom(\cM,\cN)$ consists of structure-preserving functors $F$, $G$, \ldots.
\item $\Hom(F,G)$ consists of structure-preserving natural transformations.
\end{itemize}

We unwind these definitions completely in the case of $\BrTens$ -- i.e. in the case $\cS=\LFP$.  We describe the objects, $1$-, $2$-, $3$-, and $4$-morphisms in terms of more familiar constructions in the language of braided tensor categories.  Many of the algebraic constructions which arise in this unwinding have appeared already in the study of defects in Witten--Reshetikhin--Turaev theory~\cite{Fuchs2013}, independently of the rigorous $(\infty,4)$-formulation of $\BrTens$ but clearly with this motivation in mind.  However, as there does not appear to be a source with complete details, which we require to formulate our main results, we provide them here.

\subsection{Relation to locally constant factorization algebras}\label{sec:lcfa}
The sections which follow feature a number of purely algebraic constructions, as required by Haugseng's framework.  On the other hand, it is expected that Haugseng's framework can be related to a yet-undefined framework of non-unital locally constant factorization algebras.  Hence, parallel to the rigorous extraction of the algebraic structures we need from Haugseng's general construction, we outline in Remarks \ref{rmk:BTC}, \ref{rmk:domain-wall}, \ref{rmk:domain-wall-composition} and \ref{rmk:point-defect}, and in corresponding Figures \ref{fig:BTC}, \ref{fig:domain-wall}, \ref{fig:domain-wall-composition}, and \ref{fig:point-defect} how each structure arises also in the locally constant factorization algebra framework.  This relationship is not strictly necessary to the results in the present paper, but may provide geometric intuition and motivation for our definitions.

To that end, we briefly now recall the notion of locally constant factorization algebras on topological spaces, and then its refinement to stratified spaces~\cite{Ayala2017,Costello,Ginot2015,Lurie}. A locally constant factorization algebra $\mathcal{F}$ on a topological space $X$, with coefficients in some symmetric tensor (typically, non-discrete) $1$-category $\cS$ is a pre-cosheaf, satisfying a local constancy condition, carrying a factorization structure, and finally satisfying a sheaf-like condition for gluing local sections.

Let us detail these notions one at a time. Firstly to say that $\mathcal{F}$ is a pre-cosheaf means that $\mathcal{F}$ assigns an object $\mathcal{F}(U)\in \cS$ for every open set $U\subset X$; and for every inclusion $U_1\subset U_2$, a morphism $\mathcal{F}(U_1)\to \mathcal{F}(U_2)$, compatibly with respect to composition of inclusions.  Secondly, to say that $\mathcal{F}$ is ``locally constant'' is to further require that whenever an inclusion of open sets is a retract, the induced functor is an equivalence in $\cS$.  Thirdly, factorization algebras are required to map disjoint unions of opens to tensor products in $\cS$.  This is where the name ``factorization" comes from, and it may be regarded as a ``non-commutative" multiplicative structure on $\mathcal{F}$, in that only sections which are spatially disjoint on $X$ multiply using the symmetric structure on $\cS$.  Finally, $\mathcal{F}$ must satisfy certain sheaf-like gluing conditions, but only with respect to so-called Weiss -- a.k.a. factorizing --  coverings.

More generally, given a stratified space, i.e. a filtration
$$\emptyset = X_{-1} \subset X_0 \subset X_1 \subset \cdots \subset X_n= X,$$
of $X$ by closed submanifolds $X_i$ of dimension $i$, we may consider locally constant factorization algebras with respect to the stratification.  This means the same as above, except that we only require that retracts $U_1\into U_2$ be equivalences, when $U_1$ and $U_2$ are \emph{good neighborhoods} of some closed stratum $X_i$. Here, a neighborhood is called good if for some $i$, $U\cap X_i\neq \emptyset$, and $U$ is contained in a single component of $X\backslash X_{i-1}$.

The ordinary case of locally constant factorization algebras on topological spaces is recovered as a special case when $X_i=\emptyset$ for all but the final $X_i$.

\subsection{Braided tensor categories as objects}

Following~\cite{Haugseng2017,Johnson-Freyd2017}, given a symmetric monoidal $m$-category $\cS$, there is a symmetric monoidal $(m+n)$-category of $E_n$-algebras in $\cS$.  Specializing to the case $\cS=\LFP$, and the case $n=2$, it is well known (see e.g.~\cite{Fiedorowicz1992,Fresse2017,Lurie}) that $E_2$-algebras in $\LFP$ are identified with braided tensor categories.  Hence these form the objects of $\BrTens$, as desired.

\begin{figure}\begin{center}
\begin{tikzpicture}[scale=0.75]
\begin{scope}[shift={(1,0)}]
\fill[fill=cyan!5] (-2,-2) rectangle (2,2);
\filldraw[fill=red!20, dashed,opacity=0.75] (0,0) circle (20pt);
\draw (0,0) node  {$\cA$};
\draw (2,2) node [below left] {$\RR^2$};
\draw (0,-2) node [below] {(A)};
\end{scope}
\begin{scope}[shift={(0,-5.5)}]
\fill[fill=cyan!5] (-2,-2) rectangle (2,2);
\filldraw[fill=red!20,dashed,opacity=0.75] (0,0) circle (50pt);
\filldraw[fill=red!20,dashed,opacity=0.75] (-.75,0) circle (15pt);
\filldraw[fill=red!20,dashed,opacity=0.75] (.75,0) circle (15pt);
\draw (0,-2) node [below] {(B)};
\end{scope}
\begin{scope}[shift={(6.5,0)}]
\begin{scope}
\fill[fill=cyan!5] (-2,-2) rectangle (2,2);
\filldraw[fill=red!20,dashed,opacity=0.75] (0,0) circle (50pt);
\filldraw[fill=red!20,dashed,opacity=0.75] (-.5,0) circle (30pt);
\filldraw[fill=red!20,dashed,opacity=0.75] (-1,0) circle (10pt);
\draw (-1,0) node {$1$};
\filldraw[fill=red!20,dashed,opacity=0.75] (0,0) circle (10pt);
\draw (0,0) node {$2$};
\filldraw[fill=red!20,dashed,opacity=0.75] (1,0) circle (10pt);
\draw (1,0) node {$3$};
\end{scope}

\begin{scope}[shift={(5,0)}]
\fill[fill=cyan!5] (-2,-2) rectangle (2,2);
\filldraw[fill=red!20,dashed,opacity=0.75] (0,0) circle (50pt);
\filldraw[fill=red!20,dashed,opacity=0.75] (.5,0) circle (30pt);
\filldraw[fill=red!20,dashed,opacity=0.75] (1,0) circle (10pt);
\draw (1,0) node {$3$};
\filldraw[fill=red!20,dashed,opacity=0.75] (0,0) circle (10pt);
\draw (0,0) node {$2$};
\filldraw[fill=red!20,dashed,opacity=0.75] (-1,0) circle (10pt);
\draw (-1,0) node {$1$};
\end{scope}
\draw (2.5,0) node {$\simeq$};
\draw (2.5,-2) node[below] {(C)};
\end{scope}

\begin{scope}[shift={(5.5,-5.5)}]
\begin{scope}
\fill[fill=cyan!5] (-2,-2) rectangle (2,2);
\filldraw[fill=red!20,dashed,opacity=0.75] (0,0) circle (50pt);
\filldraw[fill=red!20,dashed,opacity=0.75] (-.75,0) circle (15pt);
\draw (-.75,0) node {$1$};
\filldraw[fill=red!20,dashed,opacity=0.75] (.75,0) circle (15pt);
\draw (.75,0) node {$2$};
\end{scope}

\begin{scope}[shift={(3.5,-1)}]
\begin{scope}[scale=0.4]
\fill[fill=cyan!5] (-2,-2) rectangle (2,2);
\filldraw[fill=red!20,dashed,opacity=0.75] (0,0) circle (50pt);
\filldraw[fill=red!20,dashed,opacity=0.75] (-.75,0) circle (15pt);
\filldraw[fill=red!20,dashed,opacity=0.75] (.75,0) circle (15pt);
\draw [->, thick] (-.5,.75) arc (120:60:30pt);
\draw [->, thick] (.5,-.75) arc (300:240:30pt);
\end{scope}
\end{scope}

\begin{scope}[shift={(7,0)}]
\fill[fill=cyan!5] (-2,-2) rectangle (2,2);
\filldraw[fill=red!20,dashed,opacity=0.75] (0,0) circle (50pt);
\filldraw[fill=red!20,dashed,opacity=0.75] (-.75,0) circle (15pt);
\draw (-.75,0) node {$2$};
\filldraw[fill=red!20,dashed,opacity=0.75] (.75,0) circle (15pt);
\draw (.75,0) node {$1$};
\end{scope}

\draw [->] (2.5,0)  -- (4.5,0);
\draw (3.5,.5) node[below] {$\sim$};
\draw (3.5,-2) node[below] {(D)};
\end{scope}
\end{tikzpicture}\end{center}
\caption{\textbf{Braided tensor categories as locally constant factorization algebras on $\RR^2$.} (A) depicts a basic open set; its inclusion onto $\RR^2$ is a retract, so it is assigned the category $\cA$ canonically.  (B) depicts an embedding $\D^2\sqcup \D^2\into\D^2$; this induces the product functor $T:\cA\bt\cA\to\cA$.  (C) depicts an isotopy (with this choice of representatives, it an identity) between two composite disk embeddings $\D\sqcup\D\sqcup\D\into\D$; this induces the associator natural isomorphism $\alpha$ on $\cA$.  (D) depicts an isotopy between two disk inclusions; this induces the braiding isomorphism $\sigma$ on $\cA$.}\label{fig:BTC}
\end{figure}

\begin{remark}\label{rmk:BTC} As a special case of Lurie's theorem, the data of a braided tensor category coincides with that which defines a locally constant factorization algebra on $\RR^2$.  See Figure \ref{fig:BTC}. 
\end{remark}

\subsection{Central algebras as $1$-morphisms}\label{sec:central-algebras}
Recall from Proposition \ref{prop:AmodMonoidal} that for a braided tensor category $\cA$, the 2-category $\cA\modu$ carries a monoidal structure.  Given a pair $\cA$ and $\cB$ of braided tensor categories, we obtain similarly (by considering instead $\cA\bt\cB^{\bop}$) a monoidal structure on the 2-category of $\cA$-$\cB$ bimodule categories.  Following again~\cite{Haugseng2017,Johnson-Freyd2017}, the 1-morphisms in $\BrTens$ from $\cA$ to $\cB$ are identified with algebra objects in this monoidal 2-category.

We can capture this structure more explicitly in the notion of an $\cA$-$\cB$-central algebra, defined below.  For motivation, we recall the analogous but simpler situation in the category of vector spaces.  Given a commutative algebra $A$ over $\K$, an $A$-algebra can be defined either as an algebra object in the monoidal category $A\modu$, or equivalently as a $\K$-algebra $B$ together with an algebra morphism $A\longrightarrow Z(B)$.  This motivates the following~\cite{Drinfeld2010}:

\begin{definition-proposition}
Given a braided tensor category $\cA$, the following two notions are naturally equivalent:
\begin{itemize}
\item An $E_1$-algebra in the monoidal $2$-category of $\cA$-modules equipped with the balanced tensor product over $\cA$, and
\item A tensor category $\cC$, together with a braided tensor functor,
$$(F,J):\cA\to Z(\cC).$$
from $\cA$ to the Drinfeld center, $Z(\cC)$, of $\cC$.
\end{itemize}
We will refer to either notion as $\cA$-algebra $\cC$, or alternatively as an $\cA$-central structure on $\cC$.  \end{definition-proposition}

\begin{proof} This is a standard and straightforward verification. An algebra in $\cA\modu$ is equivalently a tensor category $\cC$ equipped with an $\cA$-balancing on the multiplication which is compatible with the associativity constraint. One checks that the functor $\cA\rightarrow \cC$ given by acting on the unit $\un_\cC$ is a tensor functor, with central structure induced by the balancing. Conversely, given a tensor functor $F:\cA \rightarrow \cC$, $\cC$ becomes an $\cA$-module via left multiplication, and a central structure on $F$ is then the same as a balancing on the multiplication of $\cC$.
\end{proof}

\noindent {\it Conventions.} (a) Given a tensor functor $F:\cA\to\cC$, an $\cA$-central structure on $\cC$ may be given by specifying half-braiding natural isomorphisms $\sigma:F(X)\ot Y \to Y \ot F(X)$.  We will call such a half-braiding natural isomorphism a central structure on $F$, and call $F$ a central functor.  (b) For a pair $\cA$, $\cB$ of braided tensor categories, we will abbreviate the phrases ``$\cA\bt\cB^{\bop}$-algebra" as ``$\cA$-$\cB$-algebra" and ``$\cA\bt\cB^{\bop}$-central structure" as ``$\cA$-$\cB$-central structure".

\begin{remark}\label{rmk:domain-wall}
In operadic terms, the triple $(\cA,\cC,F)$ is precisely the data which defines an algebra over the Swiss-Cheese operad~\cite{Voronov1999,Idrissi2017}. The data of an $\cA$-$\cB$-central algebra can also be understood in terms of locally constant factorization algebras on the stratified space $(\RR\subset\RR^2)$.  See Figure \ref{fig:domain-wall} below.
\end{remark}

\begin{example}\label{ex:centralA}
	Given a braided tensor category $\cA$, there is a canonical braided tensor functor $\cA \bt \cA^{\bop}\rightarrow Z(\cA)$ (see Proposition~\ref{prop:centerProperty2})
	Therefore $\cA$ can be seen as an $\cA$-$\cA$-algebra, a $\cA\bt\cA^{\bop}$-$\Vect$-algebra and a $\Vect$-$\cA\bt\cA^{\bop}$-algebra. In the Morita category, those are respectively the identity of $\cA$, and the evaluation and coevaluation realizing $\cA^{\bop}$ as the dual of $\cA$.
\end{example}

For later use, we record the following lemma, which is a straightforward categorification of the analogous statement for algebras over commutative algebras:
\begin{lemma}\label{lem:modmod}
Let $\cC$ be an $\cA$-algebra. Then the forgetful functor from $\cA\modu$ to $\LFP$ induces an equivalence from the category of $\cC$-modules in $\cA\modu$ to the category of $\cC$-modules in $\LFP$.
\end{lemma}
\begin{proof}
The inverse functor is given as follows.  Given a $\cC$-module $\cM$, equip it with an $\cA$-module structure using the monoidal functor $F_\cC:\cA\rightarrow \cC$. This turns $\cM$ into a $\cC$-module in $\cA\modu$, in such a way that forgetting this structure simply recovers $\cM$ as a $\cC$-module. Conversely, if $\cM$ is a $\cC$-module in $\cA\modu$, then by definition for all $a\in \cA$ and $m\in \cM$ the balanced structure on the action functor induces a natural isomorphism
	\[
		F_\cC(a)\ot m\cong a \ot m
		\]
where the left hand side uses the $\cC$-action on $\cM$, while the right hand side is the given $\cA$-action. One checks that this induces an equivalence between the given $\cA$-module structure on $\cM$ and the one induced by seeing $\cM$ as a mere category and turning it into an $\cA$-module using $F_\cC$.
\end{proof}

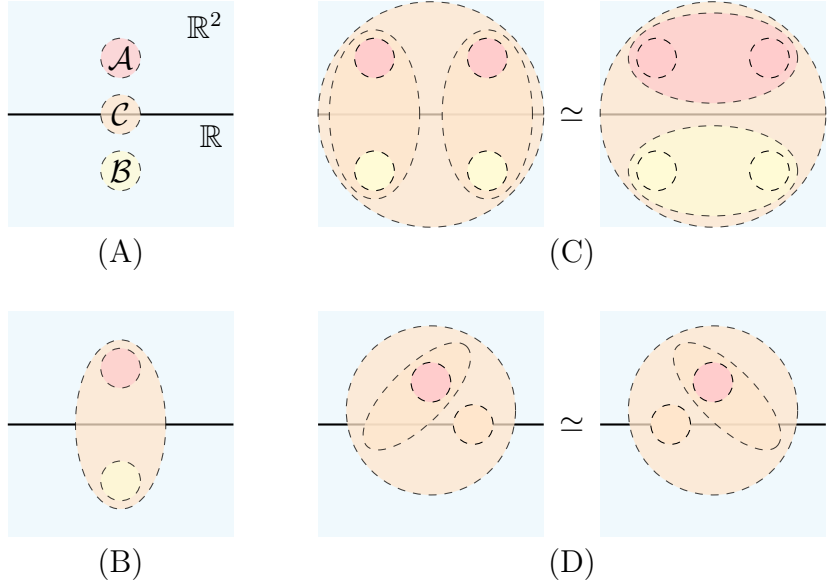
\begin{figure}\begin{center}
\begin{tikzpicture}[scale=0.75]
\begin{scope}[shift={(0,0)}]
\fill[fill=cyan!5] (-2,-2) rectangle (2,2);
\draw[thick] (-2,0) -- (2,0);
\filldraw[fill=orange!20, dashed,opacity=0.75] (0,0) circle (10pt);
\draw (0,0) node  {$\cC$};
\filldraw[fill=red!20, dashed,opacity=0.75] (0,1) circle (10pt);
\draw (0,1) node  {$\cA$};
\filldraw[fill=yellow!20, dashed,opacity=0.75] (0,-1) circle (10pt);
\draw (0,-1) node {$\cB$};
\draw (2,2) node [below left] {$\RR^2$};
\draw (2,0) node[below left] {$\RR$};
\draw (0,-2) node[below] {(A)};
\end{scope}
\begin{scope}[shift={(0,-5.5)}]
\fill[fill=cyan!5] (-2,-2) rectangle (2,2);
\draw[thick] (-2,0) -- (2,0);
\filldraw[fill=orange!20,dashed,opacity=0.75] (0,0) ellipse (.8cm and 1.5cm);
\filldraw[fill=red!20,dashed,opacity=0.75] (0,1) circle (10pt);
\filldraw[fill=yellow!20,dashed,opacity=0.75] (0,-1) circle (10pt);
\draw (0,-2) node[below] {(B)};
\end{scope}
\begin{scope}[shift={(5.5,0)}]
\begin{scope}
\fill[fill=cyan!5] (-2,-2) rectangle (2,2);
\draw[thick] (-2,0) -- (2,0);
\filldraw[fill=orange!20,dashed,opacity=0.75] (0,0) circle (2cm);
\filldraw[fill=orange!20,dashed,opacity=0.75] (-1,0) ellipse (.8cm and 1.5cm);
\filldraw[fill=orange!20,dashed,opacity=0.75] (1,0) ellipse (.8cm and 1.5cm);
\filldraw[fill=red!20,dashed] (-1,1) circle (10pt);
\filldraw[fill=yellow!20,dashed] (-1,-1) circle (10pt);

\filldraw[fill=red!20,dashed] (1,1) circle (10pt);
\filldraw[fill=yellow!20,dashed] (1,-1) circle (10pt);
\end{scope}

\begin{scope}[shift={(5,0)}]
\fill[fill=cyan!5] (-2,-2) rectangle (2,2);
\draw[thick] (-2,0) -- (2,0);
\filldraw[fill=orange!20,dashed,opacity=0.75] (0,0) circle (2cm);
\filldraw[fill=red!20,dashed,opacity=0.75] (0,1) ellipse (1.5cm and .8cm);
\filldraw[fill=yellow!20,dashed,opacity=0.75] (0,-1) ellipse (1.5cm and .8cm);

\filldraw[fill=red!20,dashed] (-1,1) circle (10pt);
\filldraw[fill=yellow!20,dashed] (-1,-1) circle (10pt);

\filldraw[fill=red!20,dashed] (1,1) circle (10pt);
\filldraw[fill=yellow!20,dashed] (1,-1) circle (10pt);
\end{scope}

\draw (2.5,0) node {$\simeq$};
\draw (2.5,-2) node[below] {(C)};
\end{scope}
\begin{scope}[shift={(5.5,-5.5)}]
\begin{scope}
\fill[fill=cyan!5] (-2,-2) rectangle (2,2);
\draw[thick] (-2,0) -- (2,0);
\filldraw[fill=orange!20,dashed,opacity=0.75] (0,.25) circle (1.5cm);
\begin{scope}[shift={(-.25,.5)}]\filldraw[fill=orange!20,dashed,opacity=0.75,rotate=45] (0,0) ellipse (1.25cm and .5cm);\end{scope}
\filldraw[fill=orange!20,dashed] (.75,0) circle (10pt);
\filldraw[fill=red!20,dashed] (0,.75) circle (10pt);
\end{scope}

\begin{scope}[shift={(5,0)}]
\fill[fill=cyan!5] (-2,-2) rectangle (2,2);
\draw[thick] (-2,0) -- (2,0);
\filldraw[fill=orange!20,dashed,opacity=0.75] (0,.25) circle (1.5cm);
\begin{scope}[shift={(.25,.5)}]\filldraw[fill=orange!20,dashed,opacity=0.75,rotate=-45] (0,0) ellipse (1.25cm and .5cm);\end{scope}
\filldraw[fill=orange!20,dashed] (-.75,0) circle (10pt);
\filldraw[fill=red!20,dashed] (0,.75) circle (10pt);

\end{scope}

\draw (2.5,0) node {$\simeq$};
\draw (2.5,-2) node[below] {(D)};
\end{scope}
\end{tikzpicture}\end{center}
\caption{\textbf{Central algebras as locally constant factorization algebras on $(\RR\subset \RR^2)$.} (A) depicts the three basic open sets $\D_\cA$, $\D_\cB$, $\D_\cC$ on $(\RR\subset\RR^2$): discs disjoint from $\RR$ are governed by the braided tensor category structure on $\cA$ and $\cB$ in each connected region, while disks intersecting $\RR$ are governed by the tensor structure of $\cC$.  (B) depicts an disks embedding $\D_\cA\sqcup\D_\cB\into \cD_\cC$; this induces a functor $F:\cA\bt\cB\to\cC$.  (C) depicts an isotopy between two disk embeddings $\D_\cA\sqcup\D_\cA\sqcup\D_\cB\sqcup\D_\cB\to\D_\cC$; this induces an isomorphism $J:F(-\ot -) \xrightarrow{\sim} F(-)\ot F(-)$, upgrading $F$ to a tensor functor.  (D) depicts an isotopy between two disk embeddings $\D_\cA\sqcup\D_\cC\to\D_\cC$; this induces a half-braiding on the image of $F$, and hence a lift $F:\cA\to Z(\cC)$.  The analogous half-braiding on $\cB$ induces a lift $F:\cB^{\bop}\to Z(\cC)$, owing to the differing orientation of the bottom region relative to $\RR$}\label{fig:domain-wall}
\end{figure}

\paragraph{Composition of $1$-morphisms}
By the preceding discussion, a $1$-morphism in $\BrTens$ from $\cB$ to $\cA$ is determined by an $\cA$-$\cB$-algebra $\cC$. Given braided tensor categories $\cA_1$ $\cA_2$, and $\cA_3$, an $\cA_1$-$\cA_2$-algebra $\cC$ and an $\cA_2$-$\cA_3$-algebra $\cD$, their composition is identified as a category with the balanced Deligne-Kelly tensor product $\cC\bt_{\cA_2}\cD$.  Let us now describe the $\cA_1$-$\cA_3$-central tensor structure on the composition.

Recall that if $\cC,\cD$ are tensor categories, then so is $\cC \bt \cD$ with multiplication given on pure tensors by
\[
	(c_1\bt d_1)\ot (c_2\bt d_2) := (c_1\ot c_2)\bt (d_1\ot b_2),
\]
and extended uniquely by cocontinuity.
\begin{definition-proposition}\label{prop:compositionCentral}
Fix a braided tensor category $\cA_2$ and central functors
	\begin{align*}
		F_\cC:\cA_2^{\bop}&\rightarrow \cC, & F_\cD:\cA_2 & \rightarrow \cD.
	\end{align*}
The composition,
	\[
\tilde T:\cC \bt \cD \bt \cC \bt \cD \longrightarrow \cC \bt \cD \longrightarrow \cC \bt_{\cA_2} \cD,
	\]
of the multiplications on $\cC$ and $\cD$ with the canonical functor (where the $\cA_2$-module structure on $\cC,\cD$ are induced by $F_\cC,F_\cD$) has a canonical balancing on the first and third tensor product. The first one is an isomorphism, for $c_i \in \cC, d_i \in \cD$, $a \in \cA$ from
\[
	\tilde T((c_1 \ot a )\bt d_1 \bt c_2 \bt d_2)=(c_1\ot a \ot c_2)\bt_{\cA_2} (d_1\ot d_2)
\]
to
\[
	\tilde T(c_1 \bt (a\ot d_1) \bt c_2 \bt d_2)=(c_1 \ot c_2)\bt_{\cA_2} (a\ot d_1\ot d_2).
\]
given by composing the central structure on $\cC$ applied on $a \ot c_2$ with the balancing on $\bt_{\cA_2}$, and the second balancing is defined similarly. Therefore, $\tilde T$ factors through a functor
	\[
		T:(\cC \bt_{\cA_2} \cD) \bt (\cC \bt_{\cA_2} \cD)  \longrightarrow \cC \bt_{\cA_2} \cD
	\]
which induces a tensor structure on $\cC \bt_{\cA_2} \cD$, such that the projection $\cC \bt\cD \rightarrow \cC \bt_{\cA_2} \cD$ is monoidal.
\end{definition-proposition}
\begin{proof}
This is proven~\cite[Theorem 6.2]{Greenough2010} in the setting of finite categories, but the proof extends to our setting as well. 
\end{proof}
\begin{definition-proposition}\label{prop:domain-wall-composition}
With the same assumptions as in Definition-Proposition \ref{prop:compositionCentral}, assume moreover that $\cC$ (resp. $\cD$) is equipped with a central functor $F$ from $\cA_1$ (resp. from $\cA_3^\bop$). Then $\cC \bt_{\cA_2}\cD$ is naturally an $\cA_1$-$\cA_3$-algebra, with central functors given by the compositions
	\[
		\cA_1 \rightarrow \cC \rightarrow \cC \bt \cD \rightarrow \cC \bt_{\cA_2} \cD
	\]
	and
	\[
		\cA_3^{\bop} \rightarrow \cD \rightarrow \cC \bt \cD \rightarrow \cC \bt_{\cA_2} \cD.
	\]
\end{definition-proposition}

Hence the composition of $1$-morphisms $\cC:\cA_2\to\cA_3$ and $\cD:\cA_1\to\cA_2$ is given by $\cC\circ\cD := \cC\bt_{\cA_2}\cD$, regarded as and $\cA_3$-$\cA_1$-central tensor category as above.

\begin{remark}\label{rmk:domain-wall-composition} In the locally constant factorization algebra framework, the balanced tensor product appearing in the composition of 2-morphisms is really a consequence of a geometric pushforward operation, together with an excision result which computes this pushforward as a balanced tensor product.  This is outlined in Figure \ref{fig:domain-wall-composition}\end{remark}

\begin{figure}\begin{center}
\begin{tikzpicture}[scale=0.75]
\begin{scope}[shift={(0,0)}]
\fill[fill=cyan!5] (-2,-3) rectangle (2,3);
\draw[thick] (-2,-1) -- (2,-1);
\draw[thick] (-2,1) -- (2,1);
\filldraw[fill=blue!20, dashed,opacity=0.75] (0,2) circle (10pt);
\draw (0,2) node {$\cA_3$};
\filldraw[fill=purple!30, dashed,opacity=0.75] (0,1) circle (10pt);
\draw (0,1) node  {$\cC$};
\filldraw[fill=red!20, dashed,opacity=0.75] (0,0) circle (10pt);
\draw (0,0) node  {$\cA_2$};
\filldraw[fill=orange!20, dashed,opacity=0.75] (0,-1) circle (10pt);
\draw (0,-1) node {$\cD$};
\filldraw[fill=yellow!20, dashed,opacity=0.75] (0,-2) circle (10pt);
\draw (0,-2) node {$\cA_1$};
\draw (0,-3) node[below] {(A)};
\end{scope}

\begin{scope}[shift={(5,0)}]
\begin{scope}
\fill[fill=cyan!5] (-2,-3) rectangle (2,3);
\draw[thick] (-2,-1) -- (2,-1);
\draw[thick] (-2,1) -- (2,1);

\end{scope}
\begin{scope}[shift={(6,0)}]
\fill[fill=cyan!5] (-2,-2) rectangle (2,2);
\draw[thick] (-2,0) -- (2,0);
\end{scope}
\draw (3,-3) node[below] {(B)};
\draw (3,0) node {$\pi$};
\draw [->] (2.25,1) -- (3.75,.25); 
\draw [->] (2.25,-1) -- (3.75,-.25); 

\end{scope}

\end{tikzpicture}\end{center}
\caption{(A) depicts the basic open sets on the stratified space $(\RR\sqcup \RR \subset\RR^2$): a locally constant factorization algebra $\cF$ is defined by labelling the basic opens as indicated.  (B) depicts a map, $\pi: (\RR\sqcup\RR\subset\RR^2) \to (\RR\subset\RR^2)$,  of stratified spaces collapsing the region between the two lines.  The composition is defined as $\cC\circ\cD := \pi_*(\cF)$.  Excision yields an equivalence $\cC\circ\cD\simeq \cC\bt_{\cA_2}\cD$, as categories, with structure maps given as in Proposition~\ref{prop:domain-wall-composition}.}\label{fig:domain-wall-composition}
\end{figure}

\subsection{Centered bimodules as $2$-morphisms}
Following again~\cite{Haugseng2017,Johnson-Freyd2017}, given two $\cA$-algebras $\cC$ and $\cD$, the 2-morphisms in $\BrTens$ are identified with $\cC$-$\cD$-bimodules internal to the monoidal $2$-category $\cA\modu$.  In this section, we introduce the notion of a \emph{centered bimodule}, which encapsulates the explicit list of functors, natural isomorphisms, and coherence conditions, which define the 2-morphisms.

\begin{definition}
Given two $\cA$-central algebras $\cC$ and $\cD$, with central functors $(F_\cC,\sigma^{\cC})$, $(F_\cD,\sigma^{\cD})$ respectively, an $\cA$-centered structure on a $\cC$-$\cD$-bimodule $\cM$ is the data of isomorphisms,
\[
	\eta_{a,m}:F_\cC(a) \ot m \cong m \ot F_\cD(a),
\]
natural in $m\in\cM$ and $a\in\cA$ and required to satisfy coherence conditions, expressed as the commutativity of the following diagrams:
\begin{equation}\label{eqn:Ax1a}
\begin{tikzpicture}
  \matrix (m) [matrix of math nodes,row sep=3em,column sep=1em,minimum width=1em]
    {
	    F_\cC(a)\ot c \ot m && c \ot m \ot F_\cD(a) \\
				& c\ot F_\cC(a) \ot m & \\};
	      \path[-stealth]
	      (m-1-1) edge node [above] {$\eta_{a,c\ot m}$} (m-1-3)
	      (m-1-1) edge node [below left] {$\sigma^{\cC}_{a,c}$}(m-2-2)
      (m-2-2) edge node [below right] {$\eta_{a,m}$}(m-1-3);
      \end{tikzpicture}
\end{equation}
\begin{equation}
\label{eqn:Ax1b}
\begin{tikzpicture}
  \matrix (m) [matrix of math nodes,row sep=3em,column sep=1em,minimum width=2em]
    {
	    m\ot d \ot F_\cD(a) && F_\cC(a) \ot m \ot d \\
				& m\ot F_\cD(a) \ot d & \\};
	      \path[-stealth]
	      (m-1-3) edge node [above] {$\eta_{a,m\ot d}$} (m-1-1)
	      (m-2-2) edge node [below left] {$\sigma^{\cD}_{a,d}$}(m-1-1)
      (m-1-3) edge node [below right] {$\eta_{a,m}$}(m-2-2);
      \end{tikzpicture}
\end{equation}

\begin{equation}
\label{eqn:Ax2}
\begin{tikzpicture}
  \matrix (m) [matrix of math nodes,row sep=3em,column sep=1em,minimum width=2em]
    {
	    F_\cC(a\ot b) \ot m && m \ot F_\cD(a \ot b) \\
    & F_\cC(b)\ot m \ot F_\cD(a) \\};
	     \path[-stealth]
	     (m-1-1) edge node [left] {$\eta_{a,b\ot m}$} (m-2-2)
	     (m-1-1) edge node [above] {$\eta_{a\ot b,m}$} (m-1-3)
	     (m-2-2) edge node [right] {$\eta_{b,m\ot a}$} (m-1-3);
      \end{tikzpicture}
\end{equation}
for $a,b\in \cA,m\in \cM, c\in\cC$ and $d\in\cD$. Here, for the sake of clarity we have omitted the explicit mention of the bimodule associator in all three diagrams.  
\end{definition}

\noindent{\bf Conventions.}  Given $\cA$-$\cB$-central algebras $\cC$ and $\cD$ and a $\cC$-$\cD$-bimodule $\cM$ we will abbreviate the phrases ``$\cA\bt\cB^{\bop}$-centered structure on $\cM$" as ``$\cA$-$\cB$-centered structure on $\cM$", and ``$\cA\bt\cB^{\bop}$-centered $\cC$-$\cD$-bimodule" as ``$\cA$-$\cB$-centered $\cC$-$\cD$-bimodule".

The following proposition identifies the notion of centered bimodule with the 2-morphisms in $\BrTens$.

\begin{proposition}
Let $\cA,\cB$ be braided tensor categories, and let $\cC,\cD$ be $\cA$-$\cB$-central tensor categories.   Then the following $2$-categories are naturally equivalent:
\begin{itemize}
\item The $2$-category of $\cA$-$\cB$-centered $\cC$-$\cD$-bimodules, and
\item The $2$-category of $\cC$-$\cD$-bimodules internal to the monoidal 2-category of $\cA$-$\cB$-bimodules.
\end{itemize}
\end{proposition}

\begin{proof}
For simplicity we treat the case of a single braided tensor category $\cA$ acting, the general case simply replaces $\cA$ by $\cA\bt\cB^{\bop}$ throughout.  We note that the first definition involves giving an isomorphism $\eta$ between two $\cA$-actions on $\cM$ -- on the left through $\cC$ and on the right through $\cD$ -- while the second definition involves a third auxiliary action of $\cA$ on $\cM$ and isomorphisms $\psi^L$ and $\psi^R$, of the left action through $\cC$ and the right action through $\cD$, respectively with the auxiliary given $\cA$-action.

Hence, starting from the second characterization, we simply compose the two isomorphisms, and check that this gives a centered structure on the bimodule.  In verifying axioms \eqref{eqn:Ax1a} and \eqref{eqn:Ax1b}, one uses that the $\cA$-bimodule structure on $\cC$ and $\cD$ invokes the central $\cA$-central structure of each algebra.  Meanwhile, axiom \eqref{eqn:Ax2} follows from the balancing axioms on $\psi^L$ and $\psi^R$.

On the other hand, starting from a centered bimodule, we can endow $\cM$ with an $\cA$-action (a left $\cA$-action through $\cC$, say), in such a way that $\psi^L$ is trivial (which we may anyway assume thanks to Lemma~\ref{lem:modmod}), and then use $\eta$ to define $\psi^R$ in the obvious way.

It is then a straightforward exercise to check that these two constructions define mutually inverse functors.
\end{proof}

\begin{example}
Let $\cC$ be an $\cA$-$\cB$-algebra. Then the $\cC$-$\cC$-bimodule $\cC$ carries the structure of an $\cA$-$\cB$-centered $\cC$-$\cC$-bimodule, induced by its central structure. This is the identity 2-morphism of $\cC$ in $\BrTens$.
\end{example}

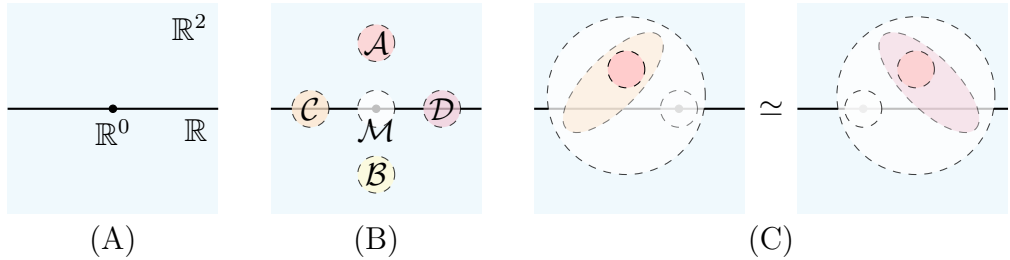
\begin{figure}[h]\begin{center}
\begin{tikzpicture}[scale=0.7]
\begin{scope}[shift={(0,0)}]
\fill[fill=cyan!5] (-2,-2) rectangle (2,2);
\draw[thick] (-2,0) -- (2,0);
\filldraw[fill=black] (0,0) circle (2pt);

\draw (2,2) node [below left] {$\RR^2$};
\draw (2,0) node[below left] {$\RR$};
\draw (0,0) node[below] {$\RR^0$};
\draw (0,-2) node[below] {(A)};
\end{scope}

\begin{scope}[shift={(5,0)}]
\fill[fill=cyan!5] (-2,-2) rectangle (2,2);
\draw[thick] (-2,0) -- (2,0);
\filldraw[fill=black] (0,0) circle (2pt);
\filldraw[fill=white!20, dashed,opacity=0.75] (-0,0) circle (10pt);
\draw (0,0) node[below] {$\cM$};
\filldraw[fill=orange!20, dashed,opacity=0.75] (-1.25,0) circle (10pt);
\draw (-1.25,0) node  {$\cC$};
\filldraw[fill=purple!20, dashed,opacity=0.75] (1.25,0) circle (10pt);
\draw (1.25,0) node  {$\cD$};
\filldraw[fill=red!20, dashed,opacity=0.75] (0,1.25) circle (10pt);
\draw (0,1.25) node  {$\cA$};
\filldraw[fill=yellow!20, dashed,opacity=0.75] (0,-1.25) circle (10pt);
\draw (0,-1.25) node {$\cB$};
\draw (0,-2) node[below] {(B)};
\end{scope}

\begin{scope}[shift={(10,0)}]
\begin{scope}
\fill[fill=cyan!5] (-2,-2) rectangle (2,2);
\draw[thick] (-2,0) -- (2,0);
\filldraw[fill=black] (.75,0) circle (2pt);
\filldraw[fill=white!20,dashed,opacity=0.75] (-.25,.25) circle (1.5cm);
\begin{scope}[shift={(-.5,.5)}]\filldraw[fill=orange!20,dashed,opacity=0.5,rotate=45] (0,0) ellipse (1.25cm and .5cm);\end{scope}
\filldraw[fill=white!20,dashed,opacity=0.5] (.75,0) circle (10pt);

\filldraw[fill=red!20,dashed] (-.25,.75) circle (10pt);
\end{scope}

\begin{scope}[shift={(5,0)}]
\fill[fill=cyan!5] (-2,-2) rectangle (2,2);
\draw[thick] (-2,0) -- (2,0);
\filldraw[fill=black] (-.75,0) circle (2pt);
\filldraw[fill=white!20,dashed,opacity=0.75] (.25,.25) circle (1.5cm);
\begin{scope}[shift={(.5,.5)}]\filldraw[fill=purple!20,dashed,opacity=0.5,rotate=-45] (0,0) ellipse (1.25cm and .5cm);\end{scope}
\filldraw[fill=white!20,dashed,opacity=0.75] (-.75,0) circle (10pt);
\filldraw[fill=red!20,dashed,opacity=0.75] (.25,.75) circle (10pt);
\end{scope}

\draw (2.5,0) node {$\simeq$};
\draw (2.5,-2) node[below] {(C)};

\end{scope}

\end{tikzpicture}\end{center}
\caption{(A) depicts the stratified space $(\RR^0\subset \RR \subset \RR^2)$.  (B) depicts the five basic open sets appearing in the stratification.  (C) depicts an isotopy between two composite disk inclusions; this induces the central structure on the bimodule $\cM$, identifying the action of $\cA$ on $\cM$ through $\cC$ and through $\cD$.  The analogous isotopy on $\cB$ gives the same structure for the induced $\cB^{\bop}$-action.}\label{fig:point-defect}
\end{figure}

\begin{remark}\label{rmk:point-defect} The data of a centered bimodule can be understood in terms of locally constant factorization alegbras on the stratified space $(\RR^0 \subset \RR^1 \subset \RR^2)$.  Note that the centered structure $\eta$ is very natural from this point of view.  See Figure \ref{fig:point-defect}, and also Example 3.2.22 of \cite{Scheimbauer2014}.\end{remark}

\paragraph{Vertical composition of $2$-morphisms}
By the preceding discussion, the $\cA$-$\cB$-centered $\cC$-$\cD$-bimodules comprise the $2$-morphisms in $\BrTens$, between the parallel $1$-morphisms $\cC, \cD:\cA\to\cB$.  Let us now turn to an explicit description of their vertical composition as 2-morphisms. Let $\cA$ be a braided tensor category, let $\cC,\cD,\cE$ be $\cA$-algebras, and let $\cM,\cN$ be $\cA$-central $\cC$-$\cD$- and $\cD$-$\cE$-bimodules. Recall that $\cM\bt_\cD \cN$ has a natural $\cC$-$\cE$-bimodule structure, where the $\cC$-action is induced by the composition \[
	\cC\bt \cM \bt \cN\xrightarrow{act} \cM \bt \cN \xrightarrow{F} \cM \bt_\cD \cN 
\] with balanced structure induced by that on $\bt_\cD$ (and likewise for the $\cE$-action).

\begin{definition-proposition}
	Denote by $F_\cC,F_\cD,F_\cE$ the central structures on $\cC,\cD,\cE$. The $\cC$-$\cE$-bimodule $\cM \bt_\cD \cN$ has a natural $\cA$-centered structure defined on pure tensors by
	\[
		(F_\cC(a)\ot m) \bt_\cD n \cong (m\ot F_\cD(a))\bt_\cD n\cong M \bt_\cD (F_\cD(a)\ot n)\cong m\bt_\cD (n \ot F_\cE(a)),
	\]
for $m\in \cM,n\in \cN, a\in \cA$, and extended uniquely by bilinearity.  Here the first map is the centered structure on $\cM$, the second is the canonical balancing on $\bt_\cD$ and the third is the centered structure on $\cN$.

\end{definition-proposition}
\begin{proof} This is a direct verification.
\end{proof}

\paragraph{Horizontal composition of 2-morphisms}
Finally, we give an explicit description of the horizontal composition of 2-morphisms. We start with the following:
\begin{proposition}
Let $\cA$ be a braided tensor category, let $\cC$, $\cD$ be $\cA$ and $\cA^{\bop}$-algebras respectively, and let $\cM$ and $\cN$ be left $\cC$ and $\cD$-modules respectively. Then the natural $\cC\bt \cD$-action on $\cM \bt \cN$ descends to an action of $\cC\bt_\cA \cD$, with its tensor structure coming from Proposition~\ref{prop:compositionCentral}, on $\cM \bt_\cA \cN$.
\end{proposition}
\begin{proof}
	The proof is the same as for Proposition~\ref{prop:compositionCentral}: the composition, 
	\[
\cC \bt\cD \bt\cM \bt \cN \rightarrow \cM \bt \cN \rightarrow \cM \bt_\cA \cN,
		\]
		of the action with the canonical functors has a canonical $\cA$-balancing on the first and third tensor product given by combining the central structure on $\cC$ and $\cD$ and the balancing on $\bt_\cA$. 
\end{proof}
\begin{corollary}\label{cor:composition2morphisms}
	Let $\cA_1,\cA_2,\cA_3$ be braided tensor categories, $\cC_1,\cD_1$ be $\cA_1$-$\cA_2$-algebras, $\cC_2,\cD_2$ be $\cA_2$-$\cA_3$-algebras and $\cM,\cN$ be centered $\cC_1$-$\cD_1$ and $\cC_2$-$\cD_2$-bimodules respectively. Then the action induced from the previous Proposition, together with the obvious centered structure, turns $\cM \bt_{\cA_2} \cN$ into a $\cA_1$-$\cA_3$-centered $\cC_1\bt_{\cA_2} \cC_2 - \cD_1 \bt_{\cA_2} \cD_2$-bimodule.
\end{corollary}
\subsection{Functors and natural transformations as $3$- and $4$-morphisms}
Using the strong (rather than lax or oplax) version of the construction in \cite{Johnson-Freyd2017}, we have the following definitions for the higher morphisms:
\begin{definition}
	A morphism $(\cM,\eta^\cM)\to(\cN,\eta^\cN)$ of $\cA$-centered $\cC$-$\cD$-bimodules is a functor 
	\[
F:\cM \longrightarrow \cN
	\]
	of $\cC$-$\cD$-bimodules, such that the following diagram commutes for $a\in \cA,m\in \cM$
\[
\begin{tikzpicture}
  \matrix (m) [matrix of math nodes,row sep=3em,column sep=5em,minimum width=1em]
    {
	    F(a\ot m)	     & F(m\ot a)\\
    a \ot F(m)&F(m)\ot a \\};
	      \path[-stealth]
	      (m-1-1) edge node [above] {$F(\eta_{a,m}^\cM)$} (m-1-2)
	      (m-2-1) edge node [below] {$\eta^\cN_{a,F(m)}$}(m-2-2)
	      (m-1-1) edge (m-2-1)
	      (m-1-2) edge (m-2-2);
      \end{tikzpicture}
\]
where the vertical arrows are the $\cC$ and $\cD$-module functor structures on $F$ respectively, and where for the sake of clarity we have suppressed the central functors.
\end{definition}

\begin{definition}
	A natural transformation of two centered bimodule functors $F,G$ is simply a natural transformation of left and right module functors.

\end{definition}
\subsection{Symmetric monoidal structure}
The 4-category $\BrTens$ has a natural symmetric monoidal structure inherited from that of $\LFP$. Let us recall here what it does on objects and morphisms in each degree.
\begin{itemize}
	\item The monoidal product of two objects $\cA,\cB$ is $\cA \bt \cB$ with its componentwise braided tensor structure.
	\item Given two braided tensor category $\cA,\cB$ and given $\cC$ and $\cD$ an $\cA$-central and a $\cB$-central tensor category, respectively, the tensor product of the central functors induces a braided tensor functor
		\[
			\cA \bt \cB \rightarrow Z(\cC)\bt Z(\cD)\simeq Z(\cC \bt \cD)
		\]
which equips $\cC\bt\cD$ with the structure of an $\cA$-$\cB$-central tensor category.  This computes the monoidal product on the 1-morphisms. 
	\item Given $\cA$-central tensor categories $\cC_1,\cD_1$, and $\cB$-central tensor categories $\cC_2,\cD_2$, an $\cA$-centered $\cC_1$-$\cD_1$-bimodule $\cM$ and a $\cB$-centered $\cC_2$-$\cD_2$-bimodule $\cN$, the Deligne--Kelly tensor product $\cM \bt \cN$ becomes a $\cC_1 \bt \cC_2$-$\cD_1\bt \cD_2$-bimodule, and the tensor product of the centered structure of $\cM$ and $\cN$ induces a central structure on $\cM \bt \cN$.  This defines the monoidal product on the 2-morphisms. 
	\item The monoidal product of functors and natural transformation is the same as in $\LFP$.
\end{itemize}

\section{The rigid Morita theories}\label{sec:rigidity}

The purpose of this section is to introduce the sub-categories $\RigidTens\subset\Tens$ and $\RigidBrTens, \BrFus\subset\BrTens$, which will turn out to be $2$-, $3$-, and $4$-dualizable categories, respectively.

\subsection{Characterizations of rigidity}
In this section, we establish some basic facts about rigidity in a tensor category, which we will use to establish dualizability.  First, we prove the equivalence between the different formulations given in the introduction:

\begin{definition-proposition}\label{def-prop:rigid}
Suppose that a tensor category $\cC$ has enough compact projectives.  Then the following conditions on $\cA$ are equivalent:
\begin{enumerate}
\item All compact projective objects of $\cC$ are left and right dualizable.
\item A generating collection of compact projective objects of $\cC$ are left and right dualizable.
\item The multiplication functor, $T: \cC\bt\cC\to\cC$ has a co-continuous right adjoint $T^R$, and the canonical lax bimodule structure on $T^R$ is strong.
\end{enumerate}
 We will say that $\cC$ is \emph{cp-rigid} if it has enough compact projectives, and if any of the above conditions is satisfied.
\end{definition-proposition}

We first start with the following, adapted from~\cite[Lemma 3.5]{Ben-Zvi2009} :
\begin{lemma}\label{lem:a-mod-dualizable} Suppose that $\cM$ and $\cN$ are $\cC$-module categories in $\LFPd$ for a tensor category $\cC$, which is generated by its compact projective objects, all of which are dualizable, and suppose that $F:\cM\to\cN$ is an $\cC$-module functor which maps compact-projective objects to compact-projectives objects, so that it admits a cocontinuous right adjoint $F^R$. Then $F^R$ has a canonical structure of a $\cC$-module functor.
\end{lemma}
\begin{proof}
Let $c=\colim c_i$ be an object in $\cC$ written as a colimit of dualizable objects, let $m\in \cM,n\in \cN$ and write $m=\colim m_i$ where the $m_i$'s are compact-projective. Since the Deligne--Kelly tensor product distributes over colimits, and since the action functor $\cC \bt \cM\rightarrow\cM$ is cocontinuous by assumption, we have
	\[
		(\colim c_i)\ot m\cong \colim (c_i \ot m).
		\]
On the other hand, if $c_i$ is dualizable we have
	\[
		\Hom_\cM(m,c_i \ot m')\cong \Hom_\cM(c_i^*\ot m,m').
		\]
		Together with the definition of an adjunction, the cocontinuity and $\cC$-linear structure of $F$ and the properties of $\Hom$ it follows that:
	\begin{align*}
		\Hom_\cM(m,F^R(\colim c_i \ot n))&\cong \Hom_\cN(\colim F(m_i),\colim c_j \ot n)\\
		&\cong \lim \Hom_\cN(F(m_i),\colim c_j \ot n)\\
		&\cong \colim \lim \Hom_\cN(F(m_i),c_j \ot n)\\
		&\cong \colim \lim \Hom_\cN(c_j^* \ot F(m_i),  n)\\
		&\cong \colim \lim \Hom_\cN(F(c_j^* \ot m_i),  n)\\
		&\cong \colim \lim \Hom_\cM(c_j^* \ot m_i, F^R( n))\\
		&\cong \colim \lim \Hom_\cM( m_i, c_i \ot F^R( n))\\
		&\cong  \Hom_\cM( m, \colim c_i \ot F^R( n))\\.
	\end{align*}
\end{proof}
We record for later use the following
\begin{corollary}\label{cor:adjointBimod}
	Under the same assumptions as in the previous Lemma, the right adjoint of a compact-preserving bimodule functor (resp. centered bimodule functor) has a canonical structure of bimodule functor (resp. of centered bimodule functor).
\end{corollary}
\begin{proof}[Proof of Definition-Proposition \ref{def-prop:rigid}]
($1\implies 2$) is clear.  The proof of ($2\implies 1$) is based on the following lemma (which uses the characterization of compact projectives given in Remark \ref{rem:projective}):

\begin{lemma}\label{lem:absoluteprojectives}
Suppose that $P=\colim P_i$ is a  small colimit of compact projective objects, and that $P$ is itself compact projective.  Then $P$ is in fact a finite direct sum of summands of the $P_i$.
\end{lemma}
\begin{proof}
First, we have an epimorphism $\oplus P_i\to\colim P_i=P$, which follows from the universal properties of each type of colimit.  Then, because $P$ is compact projective, we can split this to obtain a map $P\to \oplus_i P_i$, rendering $P$ as a summand of $\oplus_iP_i$.  Since each $P_i$ is itself compact projective, we can further split this summand as a sum of summands $P_i'$ of each $P_i$, to write $P$ as a sum of $\oplus_i P_i'$.  Finally, we can conclude that the list of nonzero $P_i'$ is finite, as otherwise $P$ would fail to be compact.
\end{proof}

\begin{remark}
The above lemma is best understood as saying that, in the sense of \cite{Brandenburg2015}, any small colimit of $2$-presentable objects (compact-projectives) is a $2$-small colimit (a direct summand of a finite direct sum) of $2$-presentable objects (compact-projectives).  There are analogous theorems about other regular cardinals.  For example, for $\aleph_0$, this says any small colimit of compact objects is a filtered colimit.
\end{remark}

Thus, having assumed a compact projective $P$ may be written a small colimit of compact projective and dualizable objects, then by the lemma, it is in fact a finite direct sum of summands of dualizable objects.  Dualizability is clearly closed under taking finite direct sums and taking direct summands, hence $P$ is dualizable.

	We may now regard $\cM=\cC\bt\cC^{\mop}$ and $\cN=\cC$ as modules for the cp-rigid tensor category $\cC^e=\cC\bt\cC^{\mop}$.  Then $T$ is a map of $\cC^e$-modules, which we claim preserves compact-projective objects. Indeed, if $x,y \in \cC$ are compact-projective, hence also dualizable by assumption, we have for $c=\colim c_i \in \cC$:
	\begin{align*}
		\Hom_\cC(x \ot y,c)&\cong \Hom_\cC(x ,(\colim c_i)\ot y^*)\\
				   & \cong \Hom_\cC(x ,\colim (c_i \ot y^*))\\
		     &\cong \colim\Hom_\cC(x , c_i \ot y^*)\\
		     &\cong \colim \Hom_\cC(x \ot y,c_i).
	\end{align*}
	
Therefore, by Lemma~\ref{lem:a-mod-dualizable}, its right adjoint $T^R$ has a canonical $\cC^e$-module structure.

Finally, for ($3\implies 1$), we can simply define,
\begin{align*}
X^* &= (\id_\cC\bt Hom(X,-))(T^R(\mathbf{1}_\cC)),\\
^*X &= (Hom(X,-)\bt\id)(T^R(\mathbf{1}_\cC)),\end{align*}
and check that they define duality functors.  For instance, we have natural isomorphisms,

\begin{align*}\Hom(\mathbf{1}_\cC, Y\ot X^*) &\cong \Hom(\mathbf{1}_\cC\bt X, (Y\bt \mathbf{1}_\cC)\ot T^R(\mathbf{1}_\cC))\\
&\cong \Hom(\mathbf{1}_\cC\bt X, T^R(Y))\\
&\cong \Hom(X, Y),
\end{align*}
from which we can extract the required evaluation and coevaluation data.

\end{proof}

\subsection{The Morita theory $\RigidTens$}	
This section is devoted to the proof of the following:
\begin{definition-proposition}
There is a well-defined symmetric monoidal 3-category $\RigidTens\subset\Tens$ whose
\begin{itemize}
\item objects are cp-rigid tensor categories,
\item 1-morphisms bimodule categories with enough compact projectives,
\item 2-morphisms are bimodule functors,
\item 3-morphisms are bimodule natural transformations.
\end{itemize}
\end{definition-proposition}

\begin{proof} Since the inclusion to $\Tens$ is full in degrees 2 and 3, we need only to show that the objects and 1-morphisms are closed under the symmetric monoidal structure, and that 1-morphisms are closed under compositions.
\paragraph{Closure under monoidal structure}
	Given two cp-rigid tensor categories, $\cC$,$\cD$, their tensor product $\cC\bt\cD$ is again a cp-rigid tensor category: it is generated under small colimits by pure tensor products of dualizable, compact and projective objects, and these are trivially again dualizable, compact and projective (see Corollary \ref{prop:purecp}).  That the monoidal product of bimodule categories with enough projectives again has enough projectives follows from the general fact that $\LFPd\subset\LFP$ is a monoidal subcategory by Proposition~\ref{prop:purecp}.

\paragraph{Closure under composition of $1$-morphisms}
Suppose that $\cC\in\RigidTens$, and that $\cM,\cN$ are right and left $\cC$-module categories, respectively, with enough compact projectives.  Then we need to show that $\cM\bt_\cC \cN$ has enough compact projectives.

For this, we recall that by construction $\cM\bt_\cC \cN$ is generated  under co-limits by the image of $\cM\bt\cN$ through the canonical functor,
$$F: \cM\bt \cN\to \cM\bt_\cC \cN.$$
Objects in the image of $F$ are called `pure tensors'.  We may rewrite this functor as
$$\cM\bt \cN \simeq (\cM\bt \cN) \bt_{\cC^e} \cC^e \xrightarrow{\id\bt m} (\cM \bt \cN)\bt_{\cC^e} \cC \simeq \cM\bt_\cC \cN.$$
By the assumption of cp-rigidity, $m$ has a co-continuous right adjoint $m^R$, hence it preserves compact projectives.  Hence it follows that $F$ preserves compact projectives as well.  Combining the two observations, the pure tensors give a collection of compact projective generators of the balanced tensor product.
\end{proof}

\subsection{The Morita theory $\RigidBrTens$}
This section is devoted to the proof of the following:
\begin{definition-proposition}
There is a well-defined symmetric monoidal 4-subcategory $\RigidBrTens\subset \BrTens$ whose:
\begin{itemize}
\item Objects are cp-rigid braided tensor categories,
\item $1$-morphisms are cp-rigid central tensor categories,
\item $2$-morphisms are centered bimodules with enough compact projectives,
\item 3- and 4-morphisms are as in $\BrTens$.
\end{itemize}
\end{definition-proposition}

\begin{proof}
As the inclusion $\RigidBrTens \subset \BrTens$ is full in degree 3 and 4, we need to check that objects, 1- and 2-morphisms are closed under the symmetric monoidal structure, and that 1- and 2-morphisms are closed under composition.

\paragraph{Closure under monoidal structure}
Again, the tensor product of two cp-rigid tensor categories is again cp-rigid, which implies the closure under tensor product of objects and 1-morphisms. Likewise, closure under the tensor product of two morphisms follows from the fact that $\LFPd$ is a monoidal sub-category of $\LFP$. 
\paragraph{Closure under composition of $1$-morphisms}
Let $\cA_1,\cA_2,\cA_3$ be cp-rigid braided tensor categories, and let $\cC$ and $\cD$ be cp-rigid $\cA_1$-$\cA_2$- and $\cA_2$-$\cA_3$-central tensor categories, respectively.  We need to show that $\cC\bt_{\cA_2} \cD$ lies in $\RigidTens$.

Since $\cA_2$ is cp-rigid, the arguments of the preceding section imply that $\cC \bt_{\cA_2} \cD$ has enough compact projectives.  The canonical functor,
	$$F: \cC\bt \cD \to \cC \bt_{\cA_2} \cD,$$
is a tensor functor, hence it maps dualizable objects to dualizable objects.  In particular, pure tensor products of dualizable, compact and projective objects will again be dualizable, compact and projective, and will be generators for $\cC\bt_{\cA_2} \cD$, as desired.

\paragraph{Closure under vertical composition of $2$-morphisms}
Let $\cA_1,\cA_2$ be cp-rigid braided tensor categories and let $\cC,\cD,\cE$ be cp-rigid $\cA_2$-$\cA_1$-central tensor categories, and let $\cM,\cN$ be $\cA_1$-$\cA_2$-centered $\cC-\cD$ and $\cD-\cE$-bimodules respectively. Then we need to show that $\cM \bt_{\cD} \cN$ has enough compact projectives.  This follows as in the case of the 1-morphisms in $\RigidTens$.
\paragraph{Closure under horizontal composition of $2$-morphisms}
The argument is the same as in the previous paragraph: it boils down to show that the balanced tensor product $\cM \bt_\cA \cN$ of two $\LFPd$-categories over a cp-rigid braided tensor category is again in $\LFPd$, which again follows as in the case of the 1-morphisms in $\RigidTens$.
\end{proof}
\subsection{The Morita theory $\BrFus$}
\begin{definition-proposition}
Over a field of characteristic zero, there is a well-defined symmetric monoidal 4-subcategory $\BrFus\subset\BrTens$ whose:
\begin{itemize}
\item Objects are braided fusion categories,
\item 1-morphisms central fusion categories.
\item 2-morphisms are centered bimodules which are the ind-completion of finite and semi-simple categories.
\item 3-morphisms are \emph{compact-preserving} functors of such.
\item 4-morphisms are natural transformations of such.
\end{itemize}
\end{definition-proposition}
\begin{proof}
	Closure for 3 and 4-morphisms is clear, so we need to check that objects, 1- and 2-morphisms are closed under the symmetric monoidal structure, and that 1- and 2- morphisms are closed under composition.  Each such claim is proved in \cite{Douglas2013,Etingof2010} in the setting of finite abelian categories.  Since the Deligne--Kelly tensor product coincides with the Deligne tensor product of finite abelian categories, and commutes with taking ind-completions, those results apply in our setting.
\end{proof}

\begin{remark}
In characteristic $p$ one needs to replace fusion with the stronger notation of ``separable" used in \cite{Douglas2013} which we briefly recall.  An algebra object $A$ in a tensor category is called separable if $A \otimes A \rightarrow A$ has a splitting as an $A$-$A$ bimodule, and a tensor category $\cC$ is called separable if the canonical algebra in $\cC \boxtimes \cC^\mop$ corresponding to the module category $\cC$ is a separable algebra.  For fusion categories, separability is equivalent to the global dimension of $\cC$ being non-zero and so is automatic in characteristic zero.  See \cite{Douglas2013} for more details.
\end{remark}

\section{Dualizability in $\Tens$ and $\BrTens$}\label{sec:dualizability}
In this section we prove Theorems \ref{thm:Tens-Informal}, \ref{thm:BrTens-Informal} and \ref{thm:BrFus-Informal}.  

We first recall what it means for an $E_1$-algebra in $\cS$ to be dualizable as a module over its enveloping algebra, and show that over such algebras, a bimodule is left and right dualizable as soon as it is dualizable as an object of $\cS$.  This step works for any $\cS$, and specializing to $\LFP$ it says that if $\cC$ is dualizable as a module over $\cC \boxtimes \cC^{\mop}$, then any bimodule category between such tensor categories with enough projectives is dualizable.

We then prove a variant of the Eilenberg-Watts theorem, which in particular establishes that module categories of the form $A\modu_\cC$, for some algebra $A\in\cC$ are dualizable, with dual ${}_\cC\modr A$.  (In \cite{Douglas2013} this is all one needs because all module categories are of this form in the finite setting, by a theorem of Ostrik~\cite{Etingof2004,Ostrik2003}.)  Finally, we show that cp-rigidity of a tensor category implies that it is dualizable over its enveloping algebra, by expressing the regular bimodule in this case in terms of an algebra object (the ``co-end", a well-known construction).  

Altogether, this shows that bimodules with enough compact projectives between cp-rigid tensor categories are dualizable.  This argument is a mix of techniques from \cite{Douglas2013} and \cite{Ben-Zvi2009,Gaitsgory2015}.  We expect that a more direct construction of the dual could be given by a single ``enriched coend" as in Remark \ref{rem:australia}, but we found the resulting category theoretic questions intimidating and so gave a more algebraic proof.

We then consider the $E_2$ setting and check that the above construction also proves dualizability for bimodules in that setting.  We conclude the section by considering dualizability of bimodule functors in the braided fusion setting.

\subsection{Dualizability over the enveloping algebra}

\begin{definition}
If $\cM$ is a left $\cC$-module in $\cS$, we say that $\cM$ is dualizable as a $\cC$-module if it is left dualizable as a $\cC$-$1$ bimodule.  Dualizability of a right module is defined similarly.
\end{definition}

\begin{definition}
If $\cC$ is an $E_1$-algebra in $\cS$, then $\cC^e = \cC \boxtimes \cC^\mop$, and $\cC$ is a left (resp. right) module over $\cC^e$ where the $\cC$ factor acts by left multiplication and $\cC^\mop$ factor acts by right multiplication.
\end{definition}

\begin{proposition}\label{prop:dualE1}
Suppose that $\cC\in \Pr$ is dualizable over $\cC^e$, and that $\cM$ is a left or right $\cC$-module, which is dualizable as an object of $\Pr$.   
Then we have a canonical equivalence,
$$\Hom_\cC(\cM,\cC)\bt_\cC \cN \simeq \Hom_\cC(\cM,\cN).$$
\end{proposition}
\begin{proof}
We treat the case that $\cM$ is a left module, as the other case is similar. First, we require the following

\begin{lemma} We have an equivalence of categories,
\begin{equation}\label{eqn:Hom-equiv}
	\Hom_{\cC^e}(\cC,\Hom(\cM,\cN)) \simeq \Hom_\cC(\cM,\cN),
\end{equation}
given by $F\mapsto F(\un_\cC)$, where $\Hom(\cM,\cN)$, -- i.e. Hom's in $\Pr$ -- is viewed as a $\cC$-bimodule by pre- and post-composition, and the $\cC^e$-module structure on the LHS uniquely determines the $\cC$-balancing on the right-hand side.
\end{lemma}
\begin{proof}
Recall that for any $\cC$-$\cC$-bimodule $\mathcal{X}$, we have a canonical identification of $\Hom_{\cC^e}(\cC,\mathcal{X})$ with the bimodule center of $\mathcal{X}$, i.e. pairs of an object of $\cX$, and a natural isomorphism between the left and right actions of $\cC$ making the obvious diagram commutes.  Applying this in the case that $\cX=\Hom(\cM,\cN)$, this latter data is equivalent to a $\cC$-linear structure on a functor, which is precisely what comprises the righthand side.  We leave the rest of the details to the reader.
\end{proof}

The result then follows from a sequence of equivalences:
\begin{align*}
	\Hom_\cC(\cM,\cC)\bt_\cC \cN&\simeq	\Hom_{\cC^e}(\cC, \Hom(\cM,\cC))\bt_\cC \cN\\
	&\simeq	\Hom_{\cC^e}(\cC,\cC^e)\bt_{\cC^e}\Hom(\cM,\cC)\bt_\cC \cN\\
	&\simeq	\Hom_{\cC^e}(\cC,\cC^e)\bt_{\cC^e}\Hom(\cM,\un_{\Pr})\bt \cC \bt_\cC \cN\\
	&\simeq	\Hom_{\cC^e}(\cC,\cC^e)\bt_{\cC^e}\Hom(\cM,\un_{\Pr}) \bt  \cN\\
	&\simeq	\Hom_{\cC^e}(\cC,\cC^e)\bt_{\cC^e}\Hom(\cM,\cN) \\
	&\simeq	\Hom_{\cC^e}(\cC,\Hom(\cM,\cN) )\\
	&\simeq	\Hom_\cC(\cM,\cN).
\end{align*}
The first and the last equivalence are \eqref{eqn:Hom-equiv}, the second follows by the assumption that $\cC$ is dualizable as a $\cC^e$-module, the third applies the dualizability of $\cM$ in $\Pr$, the fourth uses that $\cC$ is the unit for $\bt_\cC$, the fifth is again the dualizability of $\cM$ in $\Pr$, the sixth is again dualizability of $\cC$ as a $\cC^e$-module.
\end{proof}

\begin{corollary}\label{cor:separable-implies-dualizable}
Suppose that algebras $\cC,\cD$ in $\Pr$ are dualizable over their enveloping algebras, and that $\cM$ is a $\cC$-$\cD$-bimodule, which is dualizable as an object of $\Pr$.  Then $\cM$ is left and right dualizable, with
	\begin{align*}
		{}^\vee \cM&:=\Hom_\cC(\cM,\cC),&\cM^\vee&:=\Hom_\cD(\cM,\cD),
	\end{align*}
	where the $\cD$ action on ${}^\vee \cM$ is by precomposition with the right action on $\cM$, and the $\cC$ action is by postcomposition with multiplication in $\cC$ (and similarly for $\cM^{\vee}$ and $\cD$).
\end{corollary}

\begin{proof}
The argument follows the structure of \cite[\S3.2.1]{Douglas2013}. First, the canonical functor,
	\[
		\cM \bt \Hom_\cC(\cM,\cC)\rightarrow \cC,
	\]
given by evaluating the second factor on the first, has an tautological $\cD$ balanced structure, hence induces a functor
	\[
		\cM \bt_\cD \Hom_\cC(\cM,\cC)\rightarrow \cC
	\]
	which is clearly a map of $\cC$-$\cC$-bimodules.  This defines the evaluation pairing between $\cM$ and ${}^\vee\cM$.

	According to Proposition \ref{prop:dualE1}, the canonical functor
	\[
		\Hom_\cC(\cM,\cC)\bt_\cC \cM \rightarrow \Hom_\cC(\cM,\cM)
	\]
	induced by
	\[
		F \bt m \rightarrow F(-)\ot m
	\]
	is an equivalence. Hence we may define the coevaluation map as the composition
	\[
\cD \rightarrow \Hom_\cC(\cM,\cM)\simeq \Hom_\cC(\cM,\cC)\bt_\cC \cM 
	\]
where the first map is $\cD \ni d\mapsto -\ot d$ with $\cC$-linear structure given by the $\cC-\cD$-bimodule associativity constraint, which is tautologically a map of $\cD-\cD$-bimodules.

This evaluation and co-evaluation are easily shown to satisfy the zigzag relation, therefore exhibiting $\Hom_\cC(\cM,\cC)$ is a left dual for $\cM$. The proof for the right dual is similar.
\end{proof}

\begin{remark}
For simplicity, we have stated and proved Theorem \ref{prop:dualE1} and Corollary \ref{cor:separable-implies-dualizable} in the case $\cS=\Pr$, however it is clear how to adapt each statement more generally to any $\cS$: one no longer has unit objects but rather unit morphisms, and the notion of balanced tensor functor should be replaced by a morphism from the bar complex.  
\end{remark}

\subsection{Relative Eilenberg-Watts}
Let $R$ and $S$ be rings.  The Eilenberg-Watts theorem gives an equivalence between the category of co-continuous functors from $R\modu$ to $S\modu$, and the category of $S$-$R$-bimodules.  An easy consequence is that the category $R\modu$ is dualizable, with dual $\modr R$.  In this section, we establish the analogous statements for categories of co-continuous $\cC$-module functors between categories $A\modu$ and $B\modu$, where $A$ and $B$ are algebras in an arbitrary tensor category $\cC$, which we do not assume to be symmetric, or even braided.  We remind the reader that $\Hom_\cC(\cM,\cN)$ always denotes co-continuous $\cC$-module functors. 
\begin{lemma}
Let $B$ be an algebra in $\cC$, let $\cM = \mathrm{mod}_\cC$-$B$ be the left $\cC$-module category of right $B$-modules in $\cC$, and let $\cN$ be any left $\cC$-module category. Then we have an equivalence of categories,
\begin{align*}
	I_\cN:B\modu_\cN &\xrightarrow{\sim} \Hom_\cC(\cM,\cN),
\end{align*}
which maps $X\in B\modu_\cN$ to the functor $(- \ot_B X): \cM\to\cN$, equipped with the $\cC$-module structure $J$  induced by the associator,:
\[
	J_c: c \ot (-\ot_B X) \xrightarrow{\cong} (c\ot -)\ot_B X, \qquad  \textrm{for } c\in \cC,
\]
and maps a morphism $f:X\rightarrow Y$ in $B\modu_\cN$ to the $\cC$-natural transformation 
\[
\id \ot_B f:-\ot_B X \rightarrow -\ot_B Y.
\]
\end{lemma}
\begin{proof}
We show that $I_\cN$ as defined above is fully faithful and essentially surjective.  

For essential surjectivity, we note $\cM$ is generated under colimits by free modules of the form $X \ot B$, ranging over $X\in\cC$.  Since $F$ is cocontinuous, it is determined by the image of the free modules; since $F$ is $\cC$-linear it is actually determined by the image of $B$. But we have $F(B)\cong B \ot_B F(B)$ where $F(B)$ is seen as a left $B$-module via $B \ot F(B)\cong F(B\ot B) \rightarrow F(B)$, where the first isomorphism is the $\cC$-module structure on $F$. Therefore, $I_\cN$ is essentially surjective.

For faithfulness:  given $f:X\to Y$ a morphism in $B\modu_\cN$ we note that
$$I_\cN(f)(B) = (X \xrightarrow{\cong} B \ot_B X \xrightarrow{\id_B\ot_B f} B\ot_B Y \xrightarrow{\cong} Y) = f.$$
Hence, we can uniquely recover $f$ from its associated natural transformation $I_\cN(f)$.

For fullness:  given a natural transformation $\eta: I_\cN(X)\to I_\cN(Y),$ evaluation of $\eta$ on $B$ gives a morphism $f_\eta:X\to Y$, as in the proof of faithfulness; the same computation implies that $I_\cN(f_\eta)_B = \eta_B$.  Because $B$ is a $\cC$-module generator, any natural transformation between co-continuous $\cC$-module functors is determined uniquely by its value on $B\in\cM$.  Hence we can indeed conclude $I_\cN(f_\eta)=\eta$, and hence that $I_\cN$ is essentially surjective.
\end{proof}

\begin{proposition}\label{lem:mod-dualizable}The category $\cM=\modul_\cC\!-B$ is dualizable as a $\cC$-module category, with dual $\cM^\vee$ canonically identified with the category $B\modu_\cC$. 
\end{proposition}
\begin{proof}
Setting $\cN=\cC$ in the previous lemma gives an equivalence between the $\cC$-dual $\cM^\vee:=\Hom_\cC(\cM,\cC)$ and the category of left $B$-modules in $\cC$. Setting $\cN=\cM$ gives an equivalence between $\End_\cC(\cM)$ and the category of $B$-bimodules in $\cC$. Under those equivalences, the canonical functor,
\[
	\cM^\vee\boxtimes_\cC \cM \rightarrow \End_\cC(\cM)
\]
coincides with the map induced by
\[
X\boxtimes_\cC Y \mapsto X\ot Y. 
\]
Since the category of $B$-bimodules is generated by objects of this form, this map is an equivalence.
\end{proof}

\begin{remark}
It may be possible to adapt this argument to directly prove dualizability of all module categories over cp-rigid tensor categories, by replacing algebras with appropriate monads.  We will not pursue this approach here.
\end{remark}

\subsection{Cp-rigid implies dualizable over $\cC^e$}

One of the main consequences of cp-rigidity is the following:
\begin{proposition}\label{prop:rigid-implies-separable}
Any $\cC\in \RigidTens$ is dualizable over its enveloping algebra.
\end{proposition}

\begin{proof}
	This is an application of Barr-Beck monadicity theorem. We first need to show that the multiplication functor
	\[
T:\cC\bt \cC\longrightarrow \cC
		\]
has a right adjoint $T^R$ which is cocontinuous, $\cC^e$-linear and conservative.    The fact $T^R$ is cocontinuous and $\cC^e$-linear directly follows from the third characterization of cp-rigidity in Proposition~\ref{def-prop:rigid}.
	
Recall that a functor $F$ is conservative if it reflects isomorphisms, meaning that if $F(f)$ is an isomorphism then $f$ is.  To see $T^R$ is conservative, we note that $T$ is clearly essentially surjective, hence $T^R$ is faithful. Faithful functors reflect epimorphisms and monomorphisms in the above sense, but not necessarily isomorphisms.  In order to prove $F^R$ is conservative, it therefore suffices to show:

\begin{lemma}
In a category $\cC$ with enough compact projectives, every epic monomorphism is an isomorphism.\footnote{In category theory literature such a category is called `balanced'; we will avoid using this term, as it clashes with the notion of balancing of a braided monoidal category}.
\end{lemma}
\begin{proof}
Let $\Gamma \subset \cC$ denote the full sub-category of compact-projective objects. Then $\cC$ is naturally identified with the category of linear presheaves $\Lin(\Gamma^{op},\Vect)$.  Since $\Vect$ is complete and cocomplete, a natural transformation between presheaves is epic or monic if and only if it is so pointwise. Therefore, since in $\Vect$ a morphism is epic and monic if and only if it is an isomorphism, the same holds in $\cC$.\end{proof}

It then follows from Barr-Beck theorem that $\cC$ is identified with the category of modules in $\cC^e$ for the monad $T^RT$.  Following the techniques as in Theorem 4.6 of \cite{Ben-Zvi2015}, since $T^RT$ is $\cC^e$-linear we have
	\[
		T^RT(x)\cong (x \bt \un_\cC) \ot T^RT(\un_{\cC^e}).
		\]
The monad structure then turns $T^RT(\un_{\cC^e})=T^R(\un_\cC)$ into an algebra in $\cC^e$ and we have an equivalence of $\cC^e$-modules
$$\cC \simeq T^R(\un_\cC)\modu_{\cC^e}.$$
Hence $\cC$ is dualizable as a $\cC^e$-module by Lemma~\ref{lem:mod-dualizable}.
\end{proof}

\subsection{Cp-rigid tensor categories are $2$-dualizable}
Let us now apply the preceding discussion to establish the $2$-dualizability of $\RigidTens$.  
\begin{definition}
We let $h_2(B\RigidTens)$ denote homotopy 2-category of $B\RigidTens$, i.e. the $2$-category whose:
\begin{itemize}
\item Objects consist of a single object $\ast$.
\item $1$-morphisms are cp-rigid tensor categories,
\item $2$-morphisms are equivalence classes of bimodule categories with enough compact projectives
\end{itemize}
and we let $h_2(\RigidTens)$ denote the homotopy 2-category of $\RigidTens$, i.e. the 2-category whose:
\begin{itemize}
\item Objects are cp-rigid tensor categories,
\item $1$-morphisms are bimodule categories with enough compact projectives,
\item $2$-morphisms are isomorphism classes of bimodule functors.
\end{itemize}
\end{definition}

\begin{theorem}\label{thm:E1-duals}
We have:
\begin{enumerate}
\item In $h_2(B\RigidTens)$, every $1$-morphism -- i.e. every cp-rigid tensor category -- has a left and right adjoint.
\item In $h_2(\RigidTens)$, every $1$-morphism -- i.e. every bimodule category over cp-rigid categories -- has a left and a right adjoint.
\end{enumerate}
Hence, the $3$-category $\RigidTens$ is 2-dualizable.
\end{theorem}

\begin{proof}
Item (1) is well-known and applies even in $\Tens$ (and indeed any $\cS$).  Every object $\cC\in\Tens$ has a (left and right) dual given by $\cC^{\mop}$.  Note that the left and right duals are identified canonically because $\Tens$ is symmetric monoidal.  See \cite[\S3.1]{Douglas2013} for details for tensor categories, and \cite{Gwilliam2018} for the general case.

For item (2), let $\cM$ be a 1-morphism in $\RigidTens$, i.e. a $\cC$-$\cB$-bimodule category $\cM$ with enough projectives.  By Proposition \ref{prop:rigid-implies-separable}, since $\cC$ and $\cB$ are cp-rigid, they are dualizable over their enveloping algebras.  Then, by Corollary \ref{cor:separable-implies-dualizable}, the $\cC$-$\cB$-bimodule $\cM$ has both a left and a right adjoint, as required.
\end{proof}

\begin{remark} \label{rem:australia}
There should be a more direct category theoretic proof of dualizability for bimodules over cp-rigid tensor categories, which can be summarized in one sentence: it is the $\cC$-compact-projective $\cC$-enriched co-Yoneda lemma.  

In a bit more detail, we want to show that the left adjoint of ${}_\cC \cM_\cD$ is given by ${}_\cD \Hom(\cM,\cC)_\cC$.  The evaluation map is obvious and so the only issue is to write down the coevaluation.  The coevaluation is given by the following $\cC$-enriched co-end taken over the subcategory of $\cC$-compact-projectives.
$$\int^{p \in \cC-c.p.} \IHom(p, m) p.$$
Here $\IHom(p, m)$ is the internal Hom which lands in $\cC$ and $\IHom(p, m) p$ denotes the module action.
That this coend indeed gives a coevaluation then follows from the enriched co-Yoneda lemma corresponding to the enriched analogue of the compact-projective Yoneda lemma of \cite{Brandenburg2015}
$$\cF(m) = \int^{p \in \cC-c.p.} \IHom(p, m) \cF(p)$$ applied to the identity functor.
In order for this argument to work you need that there are enough $\cC$-compact-projectives, but in the cp-rigid setting any compact-projective is $\cC$-compact-projective.

The enriched co-Yoneda lemma where the enriching category is symmetric monoidal is well-known \cite{Kelly1982}.  It should be possible to modify this when the enriching category is merely monoidal, if care is taken about the difference between left-enriched and right-enriched.  The enriched Yoneda lemma in this setting appears in \cite{Hinich2016} and there is a general recipe for deriving the co-Yoneda lemma from the Yoneda lemma in any pro-arrow equipment \cite[Ex. 8.4]{Shulman2013}.  However, to our knowledge non-symmetric enriched co-Yoneda lemma has not appeared as such in the literature.  Nonetheless, there are some papers that study even more general setups.  In \cite{Street1983} the authors study enriched categories over a $2$-category, which could then be specialized to the case of $2$-categories with only one object, while in \cite[Ex. 10.2]{Garner2016} they consider $2$-categories enriched over a non-symmetric monoidal category, which can be specialized to the case where the $2$-category has only identity $2$-morphisms.  The enriched co-Yoneda lemma would then need to be modified to the $\cC$-compact-projective setting following \cite{Brandenburg2015}.  We found this task too daunting and so gave an algebraic proof that requires less heavy categorical thinking, but have included this remark because readers more adept in category theory may find the co-Yoneda construction of adjoints more natural.
\end{remark}

\subsection{Cp-rigid braided tensor categories are $3$-dualizable.}
A variation on the previous section allows us to establish the $3$-dualizability of $\RigidBrTens$.
\begin{definition}
We let $h_2(B\RigidBrTens)$ denote the 2-category whose:
\begin{itemize}
\item Objects consist of the single object $\ast$,
\item $1$-morphisms are cp-rigid braided of tensor categories,
\item $2$-morphisms are equivalence classes of central cp-rigid tensor categories,
\end{itemize}
we let $h_2(\RigidBrTens)$ denote the 2-category whose:
\begin{itemize}
\item Objects are cp-rigid braided tensor categories in $\LFPd$,
\item $1$-morphisms are central tensor categories in $\LFPd$,
\item $2$-morphisms are equivalences classes of centered bimodule categories.
\end{itemize}
and given cp-rigid braided tensor categories $\cA$ and $\cB$, we let $h_2(\Hom(\cA, \cB))$ denote the 2-category whose:
\begin{itemize}
\item Objects are $\cA$-$\cB$-central tensor categories in $\LFPd$,
\item $1$-morphisms are $\cA$-$\cB$-centered bimodule categories in $\LFPd$,
\item $2$-morphisms are isomorphism classes of centered bimodule functors.
\end{itemize}
\end{definition}

\begin{theorem}\label{thm:mainBr}
We have:
\begin{enumerate}
\item In $h_2(B\RigidBrTens)$, every 1-morphism -- i.e. every braided tensor category -- has a left and right adjoint.
\item In $h_2(\RigidBrTens)$, every 1-morphism -- i.e. every central tensor category -- has a left and a right adjoint.
\item In $h_2(\Hom(\cA,\cB))$, every 1-morphism -- i.e. every centered bimodule category -- has a left and a right adjoint.
\end{enumerate}
Hence, the $4$-category $\RigidBrTens$ is 3-dualizable.
\end{theorem}

\begin{proof}

Items (1) and (2) are general statements about $E_2$-algebras~\cite{Lurie, Calaque, Gwilliam2018}, and hold for all of $\BrTens$. 

For (1), the left and right dual to every $\cA\in\BrTens$ is $\cA^\bop$. To construct the evaluation and coevaluation maps, we observe that the identity of $\cA$, i.e. $\cA$ itself with its canonical $\cA$-$\cA$-central structure given by the functor $\cA \bt \cA^{\bop}\rightarrow Z(\cA)$ of Example~\ref{ex:centralA}, can be regarded in a tautological way as a $\cA \bt \cA^{bop}$-$\Vect$ and as a $\Vect$-$\cA^\bop \bt \cA$ central algebra using the same functor.

For (2), the left and right adjoint to any $\cA$-$\cB$-central algebra $\cC$ is $\cC^{\mop}$, with its canonical $\cB$-$\cA$-central structure induced by the equivalence $Z(\cC^{\mop})\simeq Z(\cC)^{\bop}$ from \ref{prop:centerProperty1}. Again, the evaluations and coevaluation are obtained by tautologically rewriting the centered $\cC$-bimodule $\cC$ as a centered $\cC^\mop \bt_\cA \cC$-$\cB$ bimodule and as a centered $\cA$-$\cC\bt_\cB \cC^\mop$-bimodule respectively, and similarly for the centered $\cC^\mop$-bimodule $\cC^\mop$.

Hence, it remains only to prove item (3).  So let $\cA$ and $\cB$ be cp-rigid braided tensor categories, let $\cC,\cD$ be two cp-rigid $\cA$-$\cB$-central algebras with central functors $F_\cC$ and $F_\cD$ and let $\cM$ be a centered $\cC$-$\cD$-bimodule with enough projectives.  By Theorem \ref{thm:E1-duals}, $\cM$, as a mere bimodule, has a left dual, ${}^\vee \cM =\Hom_\cC(\cM,\cC)$.  We need to equip ${}^\vee \cM$ and its evaluation and coevaluation with the structure of a centered $\cD$-$\cC$-bimodule and of maps of those respectively.  It suffices (by replacing $\cA$ by $\cA\bt\cB^{rev}$ in all arguments) to consider the case of a single braided tensor category $\cA$.

\begin{proposition}\label{prop:centered-dual}
The $\cD$-$\cC$-bimodule ${}^\vee \cM=\Hom_\cC(\cM,\cC)$ has a canonical $\cA$-centered structure given for $F\in {}^\vee\cM$ by:
\begin{multline*}
		a\ot F(-):=F(-\ot F_\cD(a))\xrightarrow{\sim} F(F_\cC(a)\ot -)\\ \xrightarrow{\sim} F_\cC(a)\ot F(-)\xrightarrow{\sim} F(- )\ot F_\cC(a)=:F(-)\ot a.
\end{multline*}
where the first map is the centered structure on $\cM$, the second is the $\cC$-module structure on $F$ and the third is the central structure on $F_\cC$.
\end{proposition}
\begin{proof}
This is direct computation.
\end{proof}

We will also need the following:
\begin{proposition}\label{prop:centeredHom}
	The $\cD$-bimodule $\Hom_\cC(\cM,\cM)$ has a canonical $\cA$-centered structure given by
\begin{multline*}
		a \ot F:= F(- \ot F_\cD(a))\xrightarrow{\sim} F(F_\cC(a) \ot -)\\ \xrightarrow{\sim} F_\cC(a)\ot F(-)\xrightarrow{\sim} F(-)\ot F_\cD(a)=: F \ot a
\end{multline*}
such that the map
	\[
		{}^\vee \cM \bt_\cC \cM \rightarrow \Hom_\cC(\cM,\cM)
	\]
given by $\varphi \bt m \mapsto \varphi(-) \otimes m$ is a morphism of centered bimodule.
\end{proposition}
\begin{proof}
This is direct computation.\end{proof}

Now the identity 2-morphism of $\cC$ is $\cC$ itself as a $\cC$-bimodule, with centered structure given by its central structure, and likewise for $\cD$. Hence, we need to show:
\begin{lemma}
The evaluation and coevaluation maps
	\begin{align*}
		\cM \bt_\cD {}^\vee \cM &\longrightarrow \cC & \cD &\longrightarrow {}^\vee \cM \bt_\cC \cM
	\end{align*}
	are morphisms of centered bimodules.
\end{lemma}
\begin{proof}
	The $\cA$-centered structure on the $\cC$-bimodule of the evaluation map is given by combining Proposition~\ref{prop:centered-dual} with Corollary~\ref{cor:composition2morphisms}. Applying the evaluation map, everything cancels except the last isomorphism of the previous Proposition, i.e. the $\cA$-central structure on $\cC$.
	recall that the coevaluation map is given by
	\[
		\cD \longrightarrow \Hom_\cC(\cM,\cM)\simeq {}^\vee \cM\bt_\cC \cM.
	\]
	The second map is an equivalence of centered bimodule by Proposition~\ref{prop:centeredHom}, so it remains to check that the first one, given by $d \mapsto - \ot d$ is a morphism of centered bimodule as well. This also follows in a straightforward way from the same Proposition, by setting $F = -\ot d$.
\end{proof}
Hence Part (2) of Theorem \ref{thm:mainBr} is proved.
\end{proof}

\subsection{Braided fusion categories are $4$-dualizable.}
Under further finiteness assumptions, we can enhance the results of the preceding section to establish $4$-dualizability in of all of $\BrFus$ over a field of characteristic zero.  (The characteristic zero assumption here is needed for this to even be a $4$-category not for dualizability, the results in this section apply to separable braided tensor categories in characteristic $p$).

\begin{definition}
Given $\cA, \cB$ braided fusion categories, and a pair $\cC,\cD$ of $\cA$-$\cB$-centered fusion categories, we let $\Hom(\cC,\cD)$ denote the 2-category whose:
\begin{itemize}
\item Objects are $\cA$-$\cB$-centered $\cC$-$\cD$-bimodules.
\item $1$-morphisms are centered bimodule compact-preserving functors.
\item $2$-morphisms are bimodule natural transformations.
\end{itemize}
\end{definition}

In addition to Theorem~\ref{thm:mainBr}, we have:

\begin{theorem}
Every $1$-morphism in $\Hom(\cC,\cD)$ -- i.e. every compact-preserving centered bimodule functor $F$ -- has a left and a right adjoint.  Hence, the $4$-category $\BrFus$ is fully dualizable.
\end{theorem}
\begin{proof}  
We follow the proof of \cite[Prop. 3.4.1]{Douglas2013} adapted to the presentable setting.  Let $F:\cM\rightarrow \cN$ be a compact-preserving centered $\cC-\cD$-bimodule functor. By definition $\cM$ and $\cN$ are ind-completions of their subcategories $\cM'$ and $\cN'$ of compact objects. Since $\cN'$ is semi-simple, every compact object in $\cN$ is in fact compact-projective, therefore $F$ preserves compact-projective objects, which implies that it arises as the ind-extension of a right exact functor
\[
F':\cM'\longrightarrow \cN'.
	\]
	Since $\cM'$ and $\cN'$ are finite and semi-simple, $F'$ has right exact left and right adjoints, which by passing to ind-completion provides compact-preserving cocontinuous left and right adjoints to $F$. It then follows from Corollary~\ref{cor:adjointBimod} that those are centered bimodules functors.
\end{proof}

\bibliographystyle{style}
\bibliography{biblio}

\end{document}